\documentclass[11pt,reqno]{amsart}

\usepackage{epsfig}
\usepackage{graphicx}
\usepackage{mathrsfs}
\usepackage{subfigure}
\usepackage{verbatim}
\usepackage{amssymb}
\usepackage{latexsym}
\usepackage{amsmath}
\usepackage{amsfonts}
\usepackage{amsbsy}
\usepackage{amscd}
\usepackage{stmaryrd}
\usepackage[dvipsnames,usenames]{color} 
\usepackage{color}
\usepackage{epstopdf}
\usepackage{hyperref}
\usepackage{multirow}
\usepackage{cite}
\usepackage{bbm}
\DeclareMathAlphabet{\mathpzc}{OT1}{pzc}{m}{it}

\makeatletter
\g@addto@macro\normalsize{%
\setlength\abovedisplayskip{6pt}
\setlength\belowdisplayskip{6pt}
\setlength\abovedisplayshortskip{6pt}
\setlength\belowdisplayshortskip{6pt}
}
\makeatother

\newcommand{\R}{{\mathbb R}}

\newtheorem{prop}{Proposition}[section]
\newtheorem{lem}[prop]{Lemma}

\newtheorem{defi}[prop]{Definition}
\newtheorem{cor}[prop]{Corollary}
\newtheorem{thm}[prop]{Theorem}

\newtheorem{rem}[prop]{Remark}

\begin{document}
\baselineskip=16pt

\title[heat equations and parabolic variational inequalities on graphs]{Semilinear heat equations and parabolic variational inequalities on graphs}
\date{\today}

\author[Y. Lin]{Yong Lin}
\author[Y. Xie]{Yuanyuan Xie$^\ast$}



\begin{abstract}
Let $G=(V,E)$ be a locally finite connected weighted graph, and $\Omega$ be an unbounded subset of $V$. Using Rothe's method, we study the existence of solutions for the semilinear heat equation $\partial_tu+|u|^{p-1}\cdot u=\Delta u~(p\ge1)$ and the parabolic variational inequality
\begin{eqnarray*}
\int_{\Omega^\circ} \partial_tu\cdot(v-u)\,d\mu\ge \int_{\Omega^\circ}(\Delta u+f)\cdot(v-u)\,d\mu
\qquad\mbox{for any }v\in \mathcal{H},
\end{eqnarray*}
where $\mathcal{H}=\{u\in W^{1,2}(V):u=0\mbox{ on }V\backslash\Omega^\circ\}$.
\end{abstract}
\maketitle

\section{Introduction}\label{S:intr}
\setcounter{equation}{0}

There is extensive literature on the study of nonlinear parabolic equation $\partial_t u+F(u)=\Delta u$.
We first introduce some results in European space. For $F(u)=|u|^{p-1}\cdot u$ and $p>1$, Gmira and Veron \cite{Gmira-Veron_1984} studied the asymptotic behaviour of solution in $\R^d$. On bounded domain $\Omega\subseteq\R^d$, Friedman and Herrero \cite{Friedman-Herrero_1987} studied extinction properties of the parabolic equation if $F(u)=\lambda u^q$ for $\lambda>0$. For $F(u)=-|u|^{\gamma-1}\cdot u$ and $\gamma>1$, Ball \cite{Ball_1977} considered the initial boundary value problem of the corresponding equation on bounded open subset of $\R^d$, and established blow-up of solution. In $\R^d$, Fujita \cite{Fujita_1966}, \cite{Fujita_1969} considered the parabolic equation when $F(u)=-u^{1+\alpha}$; he also studied the existence or nonexistence of global solutions and the blow-up of solution.

Now, we state some study about the nonlinear parabolic equation on graphs.
In 2011, Chung et al. \cite{Chung-Lee-Chung_2011} proved that the solutions of $\partial_tu+u^q=\Delta u$ on networks become extinct in finite time for $0<q<1$ and are positive for $q\ge1$. Recently, Lin and Wu \cite{Lin-Wu_2017} studied
$\partial_tu=\Delta u+u^{1+\alpha}$; using heat kernel estimates, they gave the existence and nonexistence of global solutions.
In 2018, Wu \cite{Wu_2018} proved that the semilinear heat equation $\partial_tu=\Delta_\eta u+u^{1+\alpha}$ on graph has no nonnegative global solutions under certain conditions. Also, Lin and Wu \cite{Lin-Wu_2018} studied blow-up problems for the nonlinear parabolic equation $\partial_tu=\Delta u+f(u)$ on graph.

In recent years, there is some study of other partial differential equations on graphs.
Grigoryan et al. \cite{Grigoryan-Lin-Yang_2016_1, Grigoryan-Lin-Yang_2016_2, Grigoryan-Lin-Yang_2017}, using variational method, obtained some existence results of solutions of various equations on graphs, such as Yamabe type equation, Kazdan-Warner equation and some nonlinear elliptic equations. Using Rothe's method, Lin and Xie \cite{Lin-Xie_linear_wave, Lin-Xie_nonlinear_wave} studied the existence and uniqueness of solutions of the linear and nonlinear wave equations on finite graphs.

Basing on the minimalization of the corresponding functional, variational method is commonly used to the solution of elliptic equation. However, it is not used in parabolic and hyperbolic problems, the reason is that functional with similar properties do not exist.
In 1930, Rothe originally introduced Rothe's method to study the parabolic equation. Later, this method was used to study
parabolic and hyperbolic equations (e.g., \cite{Rektorys_1971, Kacur_1984}).

A graph is an ordered pair $(V, E)$ with $V$ and $E$ being a set of vertices and edges, respectively. For any $x,y\in V$, we call $x$ a \textit{neighbour} of $y$ if $xy\in E$, and denote $x\sim y$. The \textit{degree} of $x$ is the number of its neighbours, that is,
\begin{equation*}
{\rm deg}(x):=\#\{y\in V: x\sim y\}.
\end{equation*}
A graph is called \textit{locally finite} if ${\rm deg}(x)$ is finite for any $x\in V$. Denote $M_d:=\sup_{x\in V}{\rm deg}(x)$, it can be infinite for a locally finite graph.

Let $\mu: V\to (0,\infty)$ be the vertex measure. Also, let $\omega: E\to (0,\infty)$ be the edge weight function satisfying
$$\omega_{xy}> 0, \quad\mbox{ and }\quad\omega_{xy}=\omega_{yx}\quad\mbox{ for any }xy\in E.$$
Define
\begin{eqnarray*}\label{e:D_mu}
D_\mu:=\sup\Big\{\frac{1}{\mu(x)}\sum_{y\sim x}\omega_{xy}:x\in V\Big\}.
\end{eqnarray*}
In this paper, the quadruple $G=(V,E,\mu,\omega)$ will be referred as a locally finite connected weighted graph without isolated points. Also, the graph $G$ we considered satisfies
\begin{eqnarray}\label{e:mu_0_and_D_mu}
\mu_0:=\inf_{x\in V}\mu(x)>0,\quad M_d<\infty,\quad\mbox{and}\qquad D_\mu<\infty.
\end{eqnarray}

For any nonempty subset $\Omega$ of $V$, define the boundary of $\Omega$ by
\begin{eqnarray*}\label{e:boundary_and_interior_Omega}
\partial\Omega:=\{x\in\Omega: \mbox{ there exists }y\in V\setminus\Omega\mbox{ that is a neighbour of }x\},
\end{eqnarray*}
and let $\Omega^\circ:=\Omega\backslash\partial\Omega$.

In this paper, we mainly study two problems. The first one we considered is the following semilinear heat equation
\begin{equation}\label{e:heat_eq_on_graph}
\left\{
\begin{aligned}
&\partial_tu+|u|^{p-1}\cdot u=\Delta u,\quad&&(x,t)\in\Omega^\circ\times[0,\infty),\\
&u(x,0)=h(x),\quad &&x\in\Omega^\circ,\\
&u(x,t)=0,\quad&&(x,t)\in V\setminus\Omega^\circ\times[0,\infty),
\end{aligned}
\right.
\end{equation}
where $p\ge1$ is a constant, $\Omega\subseteq V$ is unbounded, $\Delta$ is defined as in \eqref{e:defi_mu_Laplacian}, and
$h:\Omega^\circ\to\R$ is given.

Now, we define the solution of \eqref{e:heat_eq_on_graph} as follows:
\begin{defi}\label{D:solution_of_heat_eq}
A function $u: V\times[0,\infty)\to \R$ is called a {\em solution} of \eqref{e:heat_eq_on_graph} on $V\times[0,\infty)$ if $u$ is continuously differential with respect to $t$, and \eqref{e:heat_eq_on_graph} is satisfied.
\end{defi}

The following is one of our main results.

\begin{thm}\label{T:sol_of_heat_equ}
Let $G=(V, E,\mu,\omega)$ be a locally finite connected weighted graph satisfying \eqref{e:mu_0_and_D_mu}. Also, let $\Omega$ be an unbounded subset of $V$ satisfying $\Omega^\circ\neq \emptyset$. If $h\in L^2(\Omega^\circ)$, then for any $p\ge 1$, there exists a solution of \eqref{e:heat_eq_on_graph}.
\end{thm}

In $\R^d$, the solution of $u_{t}+u=\Delta u$ is given by
\begin{eqnarray*}
u(x,t)=\frac{e^{-t}}{(4\pi t)^{d/2}}\int_{\R^d}e^{-(|x-y|^2/{4t})}u_0(y)\,dy,
\end{eqnarray*}
where $u_0$ is the initial value (see \cite{Gmira-Veron_1984}). If the graph is finite, the solution of
\begin{equation}\label{e:heat_eq_on_finite_graph}
\left\{
\begin{aligned}
&\partial_tu+u=\Delta u,\quad&&(x,t)\in\Omega^\circ\times[0,\infty),\\
&u(x,0)=h(x),\quad &&x\in\Omega^\circ,\\
&u(x,t)=0,\quad&&(x,t)\in \partial\Omega\times[0,\infty),
\end{aligned}
\right.
\end{equation}
can also be explicitly given. More precisely, we get

\begin{cor}\label{C:sol_of_heat_equ_p1}
Let $G=(V, E,\mu,\omega)$ be a finite connected weighted graph, and $\Omega\subseteq V$ satisfying $\Omega^\circ\neq \emptyset$.
Let $\{\varphi_j\}_{j=1}^N$ be an orthonormal basis of $W_0^{1,2}(\Omega)$ consisting of the eigenfunctions of $-\Delta$ such that $-\Delta\varphi_j=\lambda_j\varphi_j$ for $j=1,\ldots, N$, where $N=\#\Omega^\circ$. Then the solution of \eqref{e:heat_eq_on_finite_graph} is given by
$$u(x,t)=\sum\limits_{j=1}^N h_j(0)e^{-(\lambda_j+1)t}\varphi_j(x),$$
where $h_j(0)=(h, \varphi_j)_{W_0^{1,2}(\Omega)}$.
\end{cor}

In this paper, the second problem we considered is the parabolic variational inequality. In $\R^d$, this problem has
been studied by Br\'ezis and Friedman \cite{Brezis-Friedman_1976}. In Section \ref{S:para_var_ineq}, we study
\begin{eqnarray}\label{e:para_var_ineq}
\int_{\Omega^\circ} \partial_tu\cdot(v-u)\,d\mu\ge \int_{\Omega^\circ}(\Delta u+f)\cdot(v-u)\,d\mu
\qquad\mbox{ for }t\in [0,\infty)\mbox{ and any }v\in \mathcal{H},
\end{eqnarray}
with the initial and boundary conditions
\begin{equation}\label{e:initial_boundary_value_var}
u(x,0)=g(x)\quad\mbox{ for }x\in \Omega^\circ,\quad\mbox{ and }\quad u(x,t)=0\quad\mbox{ for }(x,t)\in V\setminus\Omega^\circ\times[0,\infty),
\end{equation}
where $\Omega\subseteq V$ is unbounded, $\mathcal{H}$ is defined as in \eqref{defi:mathcal_H_m_and_H_eq}, two functions $f:\Omega^\circ\times[0,\infty)\to\R$ and $g:\Omega^\circ\to\R$ are given.

Similar to Definition \ref{D:solution_of_heat_eq}, we define the solution of
\eqref{e:para_var_ineq}--\eqref{e:initial_boundary_value_var}.
\begin{defi}\label{D:solution_of_par_var}
A function $u:V\times[0,\infty)\to\R$ is said to be a {\em solution} of \eqref{e:para_var_ineq}--\eqref{e:initial_boundary_value_var} on $V\times[0,\infty)$ if $u$ is continuously differential with respect to $t$, and \eqref{e:para_var_ineq}--\eqref{e:initial_boundary_value_var} are satisfied.
\end{defi}

Assume that there exists some positive constant $c:=c(\Omega^\circ)$ such that
\begin{equation}\label{e:condition_for_f}
\|f(\cdot,t)-f(\cdot,s)\|_{L^2(\Omega^\circ)}\le c|t-s|\qquad\mbox{for any }t,s\in[0,\infty).
\end{equation}
Now we state the last main result.
\begin{thm}\label{T:var_ineq}
Let $G=(V, E,\mu, \omega)$ be a locally finite connected weighted graph satisfying \eqref{e:mu_0_and_D_mu}, and
$\Omega$ be an unbounded subset of $V$ satisfying $\Omega^\circ\neq \emptyset$. If \eqref{e:condition_for_f} and
$g\in L^2(\Omega^\circ)$ hold, then there exists a solution of \eqref{e:para_var_ineq}--\eqref{e:initial_boundary_value_var}.
\end{thm}

The paper is organized as follows. In Section \ref{S:pre}, we state Sobolev embedding theorem on graph and theorem on existence of a solution to variational inequality. Theorem \ref{T:sol_of_heat_equ} and Corollary \ref{C:sol_of_heat_equ_p1}
are proved in Section \ref{S:proof_thm_heat_equ} and \ref{S:proof_thm_heat_equ_p1}, respectively.
In Section \ref{S:para_var_ineq}, we give the proof of Theorem \ref{T:var_ineq}.

\section{Preliminaries}\label{S:pre}
\setcounter{equation}{0}

\subsection{$\mu$-Laplacian and Sobolev embedding theorem}
Let $G=(V, E,\mu,\omega)$ be a locally finite weighted graph without isolated points.
Denote $C(V):=\{v: V\to\R\}$. The \textit{$\mu$-Laplacian} $\Delta$ of $v\in C(V)$ is defined by
\begin{equation}\label{e:defi_mu_Laplacian}
\Delta v(x)=\frac{1}{\mu(x)}\sum_{y\sim x}\omega_{xy}\big(v(y)-v(x)\big).
\end{equation}
The associated gradient form $\Gamma(\cdot,\cdot)$ and the length of $\Gamma$ are defined by
\begin{equation*}\label{e:Gamma}
\Gamma(w,v)(x)=\frac{1}{2\mu(x)}\sum_{y\sim x}\omega_{xy}\big(w(y)-w(x)\big)\big(v(y)-v(x)\big)\quad\mbox{and}\quad
|\nabla v|(x)=\sqrt{\Gamma(v,v)(x)}.
\end{equation*}
Write an integral of $v\in C(V)$ as $\int_V v\,d\mu=\sum_{x\in V}\mu(x)v(x)$.

For $U\subseteq V$, if it is finite, then we obtain the following theorem, which is one of the main tools to study solutions of the parabolic equation and the parabolic variational inequality.
\begin{lem}\label{L:Green_formula_for_mu_laplacian}({\bf Green's formula for $\mu$-Laplacian}\cite{Grigoryan_2009})
Let $G$ be a locally finite weighted graph without isolated points, and let $U\subseteq V$ be non-empty and finite. Then
\begin{eqnarray*}
\int_{U}\Delta w\cdot v\,d\mu
=-\int_U\Gamma(w,v)\,d\mu+\sum_{x\in U}\sum_{y\in V\setminus U}\omega_{xy}\big(w(y)-w(x)\big) v(x)\quad\mbox{ for any }w,v\in C(V).
\end{eqnarray*}
\end{lem}

For $1\le q<\infty$, let $L^q(U)$ be a space of all $v\in C(V)$ satisfying $\int_\Omega|v|^q\,d\mu<\infty$.
Also, let
\begin{equation}
L^\infty(U)=\{v\in C(V): \sup\limits_{x\in U}|v(x)|<\infty\}.
\end{equation}
Obviously, $L^2(U)$ is a Hilbert space with the inner product
\begin{eqnarray}\label{e:defi_inner_product_L2}
(u,v):=\int_U u v\,d\mu\qquad\mbox{ for }u,v\in L^2(U).
\end{eqnarray}
Let $W^{1,2}(U)$ be a space of all $v\in C(V)$ satisfying
\begin{eqnarray}\label{e:defi_W_12}
\|v\|_{W^{1,2}(U)}=\Big(\int_U (|\nabla v|^2+|v|^2)\,d\mu \Big)^{1/2}<\infty.
\end{eqnarray}
Let $C_0(U)$ be a set of all functions $v: U\to \R$ with $v=0$ on $\partial U$.
We let $W^{1,2}_0(U)$ be the completion of $C_0(U)$ under the norm \eqref{e:defi_W_12}.
Clearly $W^{1,2}_0(U)$ is a Hilbert space with inner product
\begin{equation}\label{eq:inner_product_W}
(w,v)_{W_0^{1,2}(U)}=\int_U(\Gamma(w,v)+wv)\,d\mu\qquad\mbox{ for any }w,v\in W_0^{1,2}(U).
\end{equation}

If $U\subseteq V$ is bounded, then it is finite. It follows that the dimension of $W_0^{1,2}(U)$ is finite,
and so we have the following embedding theorem.

\begin{thm}\label{T:Sobolev_embedding_thm}\cite{Grigoryan-Lin-Yang_2016_1}
Let $G$ be a locally finite graph, and $U\subseteq V$ be a bounded domain with non-empty interior $U^\circ$.
Then $W_0^{1,2}(U)\hookrightarrow L^q(U)$ for all $1\le q\le\infty]$. In particular, there exists a constant $C:=C(U)>0$ such that for all $1\le q\le \infty$ and all $v\in W_0^{1,2}(U)$.
\begin{eqnarray}\label{eq:PI}
\|v\|_{L^q(U)}\le C\|\nabla v\|_{L^2(U)}.
\end{eqnarray}
Moreover, $W_0^{1,2}(U)$ is pre-compact, that is, a bounded sequence in $W_0^{1,2}(\Omega)$ contains a convergent subsequence.
\end{thm}

\subsection{Existence of a solution to variational inequality}
Let $H$ be a real Hilbert space with inner product $(\cdot , \cdot)_H$ and norm $\|\cdot\|_H$.
Let $\langle\cdot, \cdot\rangle$ be the pairing between $H$ and $H'$, which is the dual of $H$.

Let $a(\cdot,\cdot)$ be a real bilinear form on $H$. We say $a(\cdot,\cdot)$ is \textit{coercive} on $H$ if there exists some constant $\beta>0$ such that
\begin{equation*}
a(v,v)\ge \beta\|v\|^2_H\quad\mbox{ for }v\in H.
\end{equation*}

\begin{thm}\label{T:Kinderlehrer-Stampacchia}\cite{Kinderlehrer-Stampacchia_1980}
Let $a(\cdot,\cdot)$ be a coercive bilinear form on $H$, and $\widetilde{H}\subset H$ be a closed and convex set. Then for given
$f\in \widetilde{H}$, there exists a unique solution of
\begin{equation*}
a(u, v-u)\ge \langle f, v-u\rangle \quad \mbox{for all }v\in \widetilde{H}.
\end{equation*}
\end{thm}

\section{Semilinear heat equation}\label{S:proof_thm_heat_equ}
\setcounter{equation}{0}

In this section, we study the existence of solution to problem \eqref{e:heat_eq_on_graph}.
Let $G=(V,E,\mu,\omega)$ be a locally finite connected weighted graph, and $\Omega$ be an unbounded subset of $V$.
Since Green's formula does not hold on unbounded set, we can not directly study problem \eqref{e:heat_eq_on_graph}. Therefore, in Subsection \ref{SS:constrction_of_bdd_subsets_eq}--\ref{SS:convergence_um_eq} we divide the research into the following steps: firstly, we construct a sequence of bounded subsets $\{\Omega_m\}$ of $\Omega$ and some needed spaces; secondly, we prove that there exists a unique solution $u_m$ of \eqref{e:heat_eq_on_graph_m}; finally, we show that $u_m$ converges to a solution $u$ of \eqref{e:heat_eq_on_graph}.

\subsection{Some spaces $\mathcal{L}^2_m,\ \mathcal{L}^2,\ \mathcal{H}_m$ and $\mathcal{H}$}\label{SS:constrction_of_bdd_subsets_eq}

For given unbounded subset $\Omega$, let $\{\Omega_m\}_{m\ge1}$ be a sequence of bounded subsets satisfying
\begin{equation}\label{e:construction_of_Omega_m_eq}
\Omega_m\subseteq \Omega_{m+1}\qquad\mbox{and}\qquad\lim_{m\to\infty}\Omega_m=\Omega.
\end{equation}
For $m\ge 1$, let $\partial\Omega_m$ and $\Omega_m^\circ$ be the boundary and interior of $\Omega_m$, respectively.
Define $L^2(V)$ and the norm as follows:
\begin{equation*}
L^2(V):=\Big\{u\in C(V): \int_V|u|^2\,d\mu<\infty\Big\},\qquad\|u\|_{L^2(V)}:=\Big(\int_V|u|^2\,d\mu\Big)^{1/2}.
\end{equation*}
Consider the subspaces $\mathcal{L}^2_m$ and $\mathcal{L}^2$ of $L^2(V)$ as
\begin{eqnarray}\label{defi:mathcal_L2_m_and_L2_eq}
\mathcal{L}^2_m=\{u\in L^2(V): u=0\mbox{ on }V\backslash\Omega^\circ_m\}\qquad\mbox{and}\qquad
\mathcal{L}^2=\{u\in L^2(V):u=0\mbox{ on }V\backslash\Omega^\circ\}
\end{eqnarray}
with norms
\begin{equation*}
\|u\|_{\mathcal{L}^2_m}:=\|u\|_{L^2(\Omega_m)}
\qquad\mbox{and}\qquad\|u\|_{\mathcal{L}^2}:=\|u\|_{L^2(\Omega)},
\end{equation*}
respectively. In view of \eqref{e:construction_of_Omega_m_eq}, we get
\begin{equation*}
\mathcal{L}^2_m\subseteq\mathcal{L}^2_{m+1}\qquad\mbox{and}\qquad \lim_{m\to\infty}\mathcal{L}^2_m=\mathcal{L}^2.
\end{equation*}
From \cite[Lemma 2.1 and Corollary 2.1]{Shao_2020}, $L^2(V)$ is a Banach space and is reflexive, and so the following results hold.
\begin{lem}\label{L:mathcal_L_and_Lm_are_closed}
$\mathcal{L}^2$  is a closed subspace of $L^2(V)$, so is $\mathcal{L}^2_m$.
Moreover, $\mathcal{L}^2$ and $\mathcal{L}^2_m$ are reflexive.
\end{lem}
\begin{proof}
Since $0\in\mathcal{L}^2$, $\mathcal{L}^2$ is nonempty. It is clear that $\mathcal{L}^2$ is a subspace of $L^2(V)$. Now we show that $\mathcal{L}^2$ is closed. For any sequence $\{v_n\}\subseteq \mathcal{L}^2$, if there exists $v\in L^2(V)$ satisfying $v_n\to v$ in $L^2(V)$, then $v\in \mathcal{L}^2$. In fact,
\begin{equation*}
\int_{V\setminus \Omega^\circ} |v_n-v|^2\,d\mu \le \int_V |v_n-v|^2\,d\mu \to 0
\end{equation*}
implies that $v=0$ on $V\setminus\Omega^\circ$. This proves $\mathcal{L}^2$ is a closed subspace of $L^2(V)$. Similarly, the result holds for $\mathcal{L}^2_m$. Since the closed subspace of the reflexive space is reflexive, we get $\mathcal{L}^2$ and $\mathcal{L}^2_m$ are reflexive.
\end{proof}

\begin{rem}\label{R:mathcal_L_is_Banach}
It follows from Lemma \ref{L:mathcal_L_and_Lm_are_closed} that $\mathcal{L}^2$ and $\mathcal{L}^2_m$ are also Banach spaces.
\end{rem}

Let $W^{1,2}(V)$ be a space of all $u\in C(V)$ satisfying $\|u\|_{W^{1,2}(V)}<\infty$. 
Let $\mathcal{H}_m$ and $\mathcal{H}$ be
\begin{eqnarray}\label{defi:mathcal_H_m_and_H_eq}
\mathcal{H}_m=\big\{u\in W^{1,2}(V): u=0\mbox{ on }V\backslash\Omega^\circ_m\big\}
\mbox{ and }\mathcal{H}=\big\{u\in W^{1,2}(V): u=0\mbox{ on }V\backslash\Omega^\circ\big\}
\end{eqnarray}
with norms
\begin{equation*}
\|u\|_{\mathcal{H}_m}:=\|u\|_{W^{1,2}(\Omega_m)}\qquad\mbox{and}\qquad\|u\|_{\mathcal{H}}:=\|u\|_{W^{1,2}(\Omega)},
\end{equation*}
respectively. From \eqref{e:construction_of_Omega_m_eq}, we get
\begin{equation*}\label{e:property_Hm_and_L2m_eq}
\mathcal{H}_m\subseteq\mathcal{H}_{m+1}\quad\mbox{and}\quad \lim_{m\to\infty}\mathcal{H}_m=\mathcal{H}.
\end{equation*}

It follows from \cite[Remark 2.2]{Shao_2020} that $W^{1,2}(V)$ is a Banach space and is reflexive. Also, $W^{1,2}$ is a Hilbert space with inner product
\begin{equation*}
(u,v)_{W^{1,2}(V)}:=\int_V(\Gamma(u,v)+uv)\,d\mu.
\end{equation*}
Consequently, the following results hold.
\begin{lem}\label{L:mathcal_H_and_Hm_are_closed}
$\mathcal{H}$ is a closed subspace of $W^{1,2}(V)$, so is $\mathcal{H}_m$. Moreover, $\mathcal{H}$ and $\mathcal{H}_m$ are reflexive.
\end{lem}
\begin{proof}
It is obvious that $0\in\mathcal{H}$ and $\mathcal{H}$ is a subspace of $W^{1,2}(V)$. For any sequence $\{v_n\}\subseteq\mathcal{H}$, if there exists $v\in W^{1,2}(V)$ satisfying $v_n\to v$ in $W^{1,2}(V)$, then $v\in\mathcal{H}$. In deed,
\begin{equation*}
\int_{V\setminus \Omega^\circ} |v_n-v|^2\,d\mu \le\int_V (|\nabla(v_n-v)|^2+|v_n-v|^2)\,d\mu \to 0.
\end{equation*}
Hence $\mathcal{H}$ is a closed subspace of $W^{1,2}(V)$.
Similarly, we obtain $\mathcal{H}_m$ is a closed subspace of $W^{1,2}(V)$. Further, $\mathcal{H}$ and $\mathcal{H}_m$ are reflexive.
\end{proof}
It is obvious that $\mathcal{H}_m$ and $\mathcal{H}$ are Hilbert spaces with the inner products
\begin{equation*}
(u, v)_{\mathcal{H}_m}=\int_{\Omega_m}\big(\Gamma(u,v)+uv\big)\,d\mu\quad\mbox{and}
\quad(u, v)_{\mathcal{H}}=\int_{\Omega}\big(\Gamma(u,v)+uv\big)\,d\mu,
\end{equation*}
respectively.

\subsection{Semilinear heat equation on $\Omega_m$.}\label{SS:initial_boundary_value_pro_eq}

In this subsection, we consider the following problem
\begin{equation}\label{e:heat_eq_on_graph_m}
\left\{
\begin{aligned}
&\partial_tu+|u|^{p-1}\cdot u=\Delta u,\quad&&(x,t)\in\Omega_m^\circ\times[0,\infty),\\
&u(x,0)=h_m(x),\quad &&x\in\Omega_m^\circ,\\
&u(x,t)=0,\quad&&(x,t)\in V\setminus\Omega_m^\circ\times[0,\infty),
\end{aligned}
\right.
\end{equation}
where $h_m=h|_{\Omega_m^\circ}$.
For $p, m\ge1$, using Rothe's method, we show that \eqref{e:heat_eq_on_graph_m} has a unique solution.

For any $T>0$, let $\{t_i\}_{i=0}^n$ be a equidistant partition of $[0,T]$ satisfying
$$t_0=0,\ t_n=T,\ \mbox{ and } t_i=i\ell\quad\mbox{ for }i\in \Lambda:=\{1,\ldots, n\}.$$
Define
\begin{equation}\label{e:defi_u_n0}
u^m_{n,0}(x)=
\left\{
\begin{aligned}
&h_m(x), \quad&&\mbox{ if }x\in\Omega_m^\circ,\\
&0,\qquad&&\mbox{ if }x\in V\setminus\Omega_m^\circ.
\end{aligned}
\right.
\end{equation}
Consider the functional $\mathcal{F}_1^m: \mathcal{H}_m\to\R$ as
\begin{eqnarray*}
\mathcal{F}^m_1(u)=\frac{1}{\ell}\int_{\Omega_m^\circ}|u|^2\,d\mu-\frac{2}{\ell}\int_{\Omega_m^\circ} u_{n,0}\cdot u\,d\mu
+\frac{2}{p+1}\int_{\Omega_m^\circ}|u|^{p+1}\,d\mu+\int_{\Omega_m}|\nabla u|^2\,d\mu.
\end{eqnarray*}
If no confusion is possible, we denote $u^m_{n,0}$ and $\mathcal{F}_1^m$ simply by $u_{n,0}$ and $\mathcal{F}_1$, respectively.
We get the following Lemma.

\begin{lem}\label{L:minimum_of_functional_eq}
$\mathcal{F}_1$ attains its minimum $u_{n,1}^m$ (or simply $u_{n,1}$) in $\mathcal{H}_m$. Also, $u_{n,1}$ is the unique solution of
\begin{eqnarray}\label{e:BVP_1}
(u-u_{n,0})/\ell+|u|^{p-1}\cdot u=\Delta u \qquad\mbox{on }\Omega_m^\circ.
\end{eqnarray}
\end{lem}

\begin{proof}
{\bf Part 1} We prove $\mathcal{F}_1$ attains its minimum in $\mathcal{H}_m$.
Obviously, H\"older inequality implies that for any $u\in \mathcal{H}_m$,
\begin{eqnarray*}
\mathcal{F}_1(u)
\ge\frac{2}{p+1}\int_{\Omega^\circ_m}|u|^{p+1}\,d\mu+\int_{\Omega_m}|\nabla u|^2\,d\mu
-\frac{1}{\ell}\int_{\Omega^\circ_m} |u_{n,0}|^2\,d\mu
\ge -\frac{1}{\ell}\|h_m\|^2_{L^2(\Omega_m^\circ)}.
\end{eqnarray*}
Thus $\mathcal{F}_1$ has a lower bound on $\mathcal{H}_m$, and then $\inf_{u\in \mathcal{H}_m}\mathcal{F}_1$ is finite. 
Taking a sequence of functions $\{u_k\}\subseteq \mathcal{H}_m$ such that $\mathcal{F}_1(u_k)\to b_1$ with $b_1=\inf_{u\in\mathcal{H}_m} \mathcal{F}_1(u)$. This leads to $|\mathcal{F}_1(u_k)-b_1|<\epsilon_1$ for some $\epsilon_1>0$. As a result,
\begin{eqnarray*}
\int_{\Omega_m}|\nabla u_k|^2\,d\mu
\le b_1+\epsilon_1+\frac{1}{\ell}\int_{\Omega^\circ_m}|u_{n,0}|^2\,d\mu.
\end{eqnarray*}
Combining this with Theorem \ref{T:Sobolev_embedding_thm}, we get $u_k$ is bounded in $\mathcal{H}_m$. Moreover, there exist a function $u_{n,1}\in \mathcal{H}_m$
and a subsequence of $\{u_{k}\}$, denoted by itself, such that $u_{k}\to u_{n,1}$ in $\mathcal{H}_m$. So
$\|u_{k}\|_{W^{1,2}(\Omega_m)}\to \|u_{n,1}\|_{W^{1,2}(\Omega_m)}$. In view of
\begin{eqnarray*}
\big|\|u_{k}\|_{L^2(\Omega_m)}-\|u_{n,1}\|_{L^2(\Omega_m)}\big|
\le \|u_{k}-u_{n,1}\|_{L^2(\Omega_m)}
\le \|u_{k}-u_{n,1}\|_{W^{1,2}(\Omega_m)},
\end{eqnarray*}
we get
\begin{eqnarray}\label{e:u_nk_to_u_infty_in_L2}
\|u_{k}\|_{L^2(\Omega_m)}\to \|u_{n,1}\|_{L^2(\Omega_m)}\quad\mbox{and}\quad
\|\nabla u_{k}\|_{L^2(\Omega_m)}\to \|\nabla u_{n,1}\|_{L^2(\Omega_m)}.
\end{eqnarray}
Furthermore, $u_{k}\to u_{n,1}$ on $\Omega_m$. Hence
\begin{eqnarray*}
\mathcal{F}_1(u_{n,1})=\lim_{k\to\infty}\mathcal{F}_1(u_{k})=b_1.
\end{eqnarray*}
That is, $\mathcal{F}_1$ attains its minimum $u_{n,1}$ in $\mathcal{H}_m$.

{\bf Part 2.} We gave that $u_{n,1}$ is the unique solution of \eqref{e:BVP_1}.
For any $\phi\in\mathcal{H}_m$,
\begin{eqnarray*}\label{e:Eular_eq}
\begin{aligned}
0=&\frac{d}{d\xi}\Big|_{\xi=0}\mathcal{F}_1(u_{n,1}+\xi\varphi)\\
=&\frac{2}{\ell}\int_{\Omega^\circ_m}(u_{n,1}-u_{n,0})\cdot \phi\,d\mu
+2\int_{\Omega^\circ_m}|u_{n,1}|^{p-1}\cdot u_{n,1}\cdot\phi\,d\mu
-2\int_{\Omega_m^\circ}\Delta u_{n,1}\cdot \phi\,d\mu,
\end{aligned}
\end{eqnarray*}
where we use the fact that for any $v_1,v_2\in\mathcal{H}_m$,
\begin{equation}\label{e:green_formula_eq}
-\int_{\Omega^\circ_m}\Delta v_1\cdot v_2\,d\mu=\int_{\Omega_m}\Gamma(v_1,v_2)\,d\mu.
\end{equation}
The derivation is as follows. Using Lemma \ref{L:Green_formula_for_mu_laplacian}, we get for any $v_1,v_2\in\mathcal{H}_m$,
\begin{eqnarray*}
\begin{aligned}
&-2\int_{\Omega_m^\circ}\Delta v_1\cdot v_2\,d\mu\\
=&\sum_{x\in\Omega_m^\circ}\Big[\sum_{y\sim x, y\in \Omega^\circ_m}\omega_{xy}\big(v_1(y)-v_1(x)\big)\big(v_2(y)-v_2(x)\big)
+\sum_{y\sim x, y\in \partial\Omega_m}\omega_{xy}\big(v_1(y)-v_1(x)\big)\big(v_2(y)-v_2(x)\big)\Big]\\
&+\sum_{y\in\Omega^\circ_m}\sum_{x\sim y, x\in \partial\Omega_m}\omega_{yx}\big(v_1(x)-v_1(y)\big)\big(v_2(x)-v_2(y)\big)\\
=&\sum_{x\in\Omega_m^\circ}\sum_{y\sim x, y\in \Omega_m}\omega_{xy}\big(v_1(y)-v_1(x)\big)\big(v_2(y)-v_2(x)\big)
+\sum_{x\in\partial\Omega_m}\sum_{y\sim x, y\in \Omega^\circ_m}\omega_{xy}\big(v_1(y)-v_1(x)\big)\big(v_2(y)-v_2(x)\big)\\
=&\sum_{x\in\Omega_m}\sum_{y\sim x, y\in \Omega_m}\omega_{xy}\big(v_1(y)-v_1(x)\big)\big(v_2(y)-v_2(x)\big)\\
=&2\int_{\Omega_m}\Gamma(v_1,v_2)\,d\mu.
\end{aligned}
\end{eqnarray*}

If $\widetilde{u}\in \mathcal{H}_m$ is another solution of \eqref{e:BVP_1}, then
\begin{eqnarray}\label{e:equation_on_solution_u1}
1/\ell\cdot(u_{n,1}-\widetilde{u})-\Delta(u_{n,1}-\widetilde{u})
+\big(|u_{n,1}|^{p-1}\cdot u_{n,1}-|\widetilde{u}|^{p-1}\cdot \widetilde{u}\big)=0\qquad \mbox{on }\Omega_m^\circ.
\end{eqnarray}
According to \cite{Grigoryan_2009}, $-\Delta$ is a positive definite Laplacian. For any $x_0\in \Omega_m^\circ$,
if $(u_{n,1}-\widetilde{u})(x_0)\ge0$, then
\begin{equation*}
-\Delta(u_{n,1}-\widetilde{u})(x_0)\ge 0 \quad\mbox{and}\quad
\big(|u_{n,1}|^{p-1}\cdot u_{n,1}-|\widetilde{u}|^{p-1}\cdot\widetilde{u}\big)(x_0)\ge 0,
\end{equation*}
which, together with \eqref{e:equation_on_solution_u1}, yield $u_{n,1}(x_0)=\widetilde{u}(x_0)$.
Since $x_0$ is arbitrary, we obtain $u_{n,1}=\widetilde{u}$ on $\Omega^\circ$, which completes the proof.
\end{proof}

Successively, for $i=2,\ldots,n$, consider the functionals $\mathcal{F}^m_i (\mbox{or simply }\mathcal{F}_i) : \mathcal{H}_m\to \R$:
\begin{eqnarray*}
\mathcal{F}_i(u)=\frac{1}{\ell}\int_{\Omega^\circ_m}|u|^2\,d\mu-\frac{2}{\ell}\int_{\Omega^\circ_m} u_{n,i-1}\cdot u\,d\mu
+\frac{2}{p+1}\int_{\Omega^\circ_m}|u|^{p+1}\,d\mu+\int_{\Omega_m}|\nabla u|^2\,d\mu.
\end{eqnarray*}
Similarly, $\mathcal{F}_i$ attains its minimum $u^m_{n,i}$ (or simply $u_{n,i}$) in $\mathcal{H}_m$
that is the unique solution of
\begin{eqnarray}\label{e:BVP_i}
(u-u_{n,i-1})/\ell+ u\cdot |u|^{p-1}=\Delta u\qquad\mbox{on }\Omega^\circ_m.
\end{eqnarray}
Let $u_{n,i}(x)$ be the approximation of a solution $u_m(x,,t)$ of \eqref{e:heat_eq_on_graph_m} at $t=t_i$.
Define
\begin{equation*}\label{e:defi_delta_ui}
\delta u^m_{n,i}(x)=(u_{n,i}(x)-u_{n,i-1}(x))/\ell\quad\mbox{ for }i\in \Lambda,
\end{equation*}
and denote $\delta u^m_{n,i}(x)$ simply by $\delta u_{n,i}(x)$.
By \eqref{e:BVP_1} and \eqref{e:BVP_i}, we obtain for $i\in \Lambda$,
\begin{equation}\label{e:delta2_ui_is_a_uni_sol}
\delta u_{n,i}+|u_{n,i}|^{p-1}\cdot u_{n,i}=\Delta u_{n,i}\qquad\mbox{on }\Omega^\circ_m.
\end{equation}
Define Rothe's sequence $\{u^{(n),m}(x,t)\}$ (or simply $\{u^{(n)}(x,t)\}$) from $[0,T]$ to $\mathcal{H}_m$ by
\begin{equation}\label{e:defi_u(n)}
u^{(n)}(x,t)=u_{n,i-1}(x)+(t-t_{i-1})\cdot\delta u_{n,i}(x)\quad t\in[t_{i-1}, t_i].
\end{equation}
Also, define auxiliary function as
\begin{equation}\label{e:defi_delta_overline_u(n)}
\overline{u}^{(n)}(x,t)=
\left\{
\begin{aligned}
&u_{n,i}(x),\qquad&&(x,t)\in\Omega_m^\circ\times(t_{i-1},t_i],\\
&h_m(x),\qquad&&(x,t)\in\Omega_m^\circ\times[-\ell,0],\\
&0,\qquad&&(x,t)\in V\setminus\Omega^\circ_m\times[-\ell,0].
\end{aligned}
\right.
\end{equation}
In the following, we give some priori estimates, which will be used to show that Rothe's function $u^{(n)}(x,t)$
converges to a solution $u_m(x,t)$ of \eqref{e:heat_eq_on_graph_m}.

\begin{lem}\label{L:est_of_u(n)_et_eq}
There exist positive constants $C_j$, $j=1,2,3,4$, such that for $i\in \Lambda$,
\begin{eqnarray}\label{e:bounded_of_delta ui_and_ui_and_delta2_ui}
\|u_{n,i}\|^2_{L^2(\Omega_m)}\le C_1,
\quad \|\nabla u_{n,i}\|^2_{L^2(\Omega_m)}\le C_2,
\quad\|u_{n,i}\|^2_{L^{2p}(\Omega_m)}\le C_3,
\quad\|\delta u_{n,i}\|^2_{L^2(\Omega_m)}\le C_4,
\end{eqnarray}
where $C_j:=C_j(\Omega_m)$ for $j=1,2,3$, and $C_4:=C_4(\Omega_m, p)$.
\end{lem}

\begin{proof}
It follows from \eqref{e:delta2_ui_is_a_uni_sol} that for any $v\in \mathcal{H}_m$,
\begin{equation*}
\int_{\Omega_m^\circ}\delta u_{n,i}\cdot v\,d\mu+\int_{\Omega_m^\circ}|u_{n,i}|^{p-1}\cdot u_{n,i}\cdot v\,d\mu
=\int_{\Omega_m^\circ}\Delta u_{n,i}\cdot v\,d\mu.
\end{equation*}
Letting $v=u_{n,i}$ in the above equality, using H\"older inequality and \eqref{e:green_formula_eq}, we get for any $i\in \Lambda$,
\begin{eqnarray*}
\begin{aligned}
\int_{\Omega^\circ_m}|u_{n,i}|^2\,d\mu
\le&\int_{\Omega^\circ_m}|u_{n,i}|^2\,d\mu+\ell\Big(\int_{\Omega^\circ_m}|u_{n,i}|^{p+1}\,d\mu
+\int_{\Omega_m}|\nabla u_{n,i}|^2\,d\mu\Big)\\
\le&\frac{1}{2}\Big(\int_{\Omega^\circ_m}|u_{n,i}|^2\,d\mu+\int_{\Omega^\circ_m}|u_{n,i-1}|^2\,d\mu\Big).
\end{aligned}
\end{eqnarray*}
This leads to for any $i\in \Lambda$,
\begin{eqnarray*}
\|u_{n,i}\|^2_{L^2(\Omega_m)}\le\|u_{n,i-1}\|^2_{L^2(\Omega_m)}\le\cdots\le\|u_{n,0}\|^2_{L^2(\Omega_m)}
=\|h_m\|^2_{L^2(\Omega_m^\circ)}=:C_1.
\end{eqnarray*}
Furthermore, $|u_{n,i}|\le \sqrt{C_1/\mu_0}$ on $\Omega_m$. Hence
\begin{eqnarray*}
\|\nabla u_{n,i}\|^2_{L^2(\Omega_m)}\le 2D_\mu \cdot|u_{n,i}(x)|\cdot \mu(\Omega_m)\le C_2.
\end{eqnarray*}
It follows from Theorem \ref{T:Sobolev_embedding_thm} that
\begin{eqnarray*}
\|u_{n,i}\|^2_{L^{2p}(\Omega_m)}\le C(\Omega_m)\cdot\|\nabla u_{n,i}\|^2_{L^2(\Omega_m)}\le C_3.
\end{eqnarray*}

In view of \eqref{e:delta2_ui_is_a_uni_sol}, using $(a-b)^2\le 2(a^2+b^2)$ and
$|\Delta u_{n,i}(x)|^2\le D_\mu|\nabla u_{n,i}(x)|^2$, we get
\begin{eqnarray*}
\|\delta u_{n,i}\|^2_{L^2(\Omega_m)}\le2D_\mu\|\nabla u_{n,i}(x)\|^2_{L^2(\Omega_m)}+2\|u_{n,i}\|^{2p}_{L^{2p}(\Omega_m)}\le C_4.
\end{eqnarray*}
This completes the proof.
\end{proof}
According to Lemma \ref{L:est_of_u(n)_et_eq}, we can derived directly the following results.
\begin{lem}\label{L:est_of_u(n)_et}
There exist constants $C_5:=C_5(\Omega_m, p)$, $C_6:=C_6(\Omega_m)$ such that for any $t\in[0,T]$,
\begin{eqnarray}\label{e:bounded_1}
\|u^{(n)}(\cdot,t)\|_{L^{2}(\Omega_m)}+\|\overline{u}^{(n)}(\cdot,t)\|_{L^{2}(\Omega_m)}
+\|\overline{u}^{(n)}(\cdot,t)\|_{L^{2p}(\Omega_m)}+\|\partial_tu^{(n)}(\cdot,t)\|_{L^2(\Omega_m)}
\le C_5,
\end{eqnarray}
\begin{eqnarray}\label{e:bounded_2}
\|u^{(n)}(\cdot,t)-\overline{u}^{(n)}(\cdot,t)\|_{L^2(\Omega_m)}\le C_6/n.
\end{eqnarray}
\end{lem}

\begin{lem}\label{L:weakly_conv}
There exist a function $u_m\in \mathcal{L}^2_m$ satisfying $\partial_t u_m\in\mathcal{L}^2_m$,
and two subsequences of $\{u^{(n)}\}$, $\{\overline{u}^{(n)}\}$, denoted by themselves, such that for any $(x,t)\in\Omega^\circ_m\times[0,T]$,
\begin{enumerate}
\item[(a)] $u^{(n)}\to u_m$ and $\overline{u}^{(n)}\to u_m$;
\item[(b)] $\partial_tu^{(n)}\to \partial_tu_m$.
\end{enumerate}
\end{lem}

\begin{proof}
(a) According to \eqref{e:bounded_1}, we get there exist two subsequences $\{u^{(n)}\}, \{\overline{u}^{(n)}\}$ and two functions $u_m, \overline{u}_m\in \mathcal{L}^2_m$ satisfying for $t\in[0,T]$,
\begin{eqnarray}\label{e:inner_conve_unk_to_u}
u^{(n)}(\cdot,t)\to u_m(\cdot,t)\quad\mbox{and}
\quad\overline{u}^{(n)}(\cdot,t)\to\overline{u}_m(\cdot,t)\quad\mbox{ on }\mathcal{L}^2(\Omega_m).
\end{eqnarray}
This leads to
\begin{eqnarray*}
u^{(n)}\to u_m\quad\mbox{and}\quad\overline{u}^{(n)}\to\overline{u}_m\quad\mbox{on }\Omega^\circ_m\times[0,T],
\end{eqnarray*}
which, together with \eqref{e:bounded_2}, yield
\begin{eqnarray*}
 \|u_m(\cdot,t)-\overline{u}_m(\cdot,t)\|^2_{L^2(\Omega_m)}
=\lim_{n\to\infty}\|u^{(n)}(\cdot,t)-\overline{u}^{(n)}(\cdot,t)\|^2_{L^2(\Omega_m)}=0.
\end{eqnarray*}
Hence $u_m=\overline{u}_m$ on $\Omega^\circ_m\times[0,T]$. This proves part (a).

(b) By \eqref{e:bounded_1}, we can extract a subsequence of $\{\partial_tu^{(n)}\}$ such that for a function $S_m\in\mathcal{L}^2_m$,
\begin{equation*}
\partial_tu^{(n)}\to S_m\qquad\mbox{on }\Omega_m^\circ\times[0,T],
\end{equation*}
where $S_m=\partial_tu_m$. In deed, for any $t\in [0,T]$,
\begin{eqnarray}\label{e:int_par_s_u_n_eq}
u^{(n)}(x,t)-h_m(x)=\int_0^t\partial_s u^{(n)}(x,s)\,ds\quad\mbox{ on }\Omega_m^\circ.
\end{eqnarray}
It follows from \eqref{e:bounded_1} that
\begin{eqnarray*}
|\partial_t u^{(n)}(x,t)|\le C_5/{\sqrt{\mu_0}}\quad\mbox{ on }\Omega_m^\circ\times[0,T].
\end{eqnarray*}
Using Dominated Convergence Theorem, we get
\begin{equation*}
\int_0^t\partial_su^{(n)}(x,s)\,ds\to \int_0^t S_m(x,s)\,ds\quad\mbox{ on }\Omega_m^\circ\times[0,T].
\end{equation*}
Taking the limits as $n\to\infty$ in \eqref{e:int_par_s_u_n_eq}, we get
\begin{eqnarray*}
u_m(x,t)-h_m(x)=\int_0^tS_m(x,s).
\end{eqnarray*}
Hence $S_m=\partial_t u_m$. Also, $u_m(x,0)=h_m(x)$ on $\Omega^\circ_m$.
\end{proof}

\begin{lem}\label{L:int_0_T_conv}
The following results hold on $\Omega_m^\circ$:
\begin{enumerate}
\item[(a)] $\int_0^T \Delta \overline{u}^{(n)}\,dt\to\int_0^T\Delta u_m\,dt$;
\item[(b)] $\int_0^T|\overline{u}^{(n)}|^{p-1}\cdot\overline{u}^{(n)}\,dt\to\int_0^T|u_m|^{p-1}\cdot u_m\,dt$.
\end{enumerate}
\end{lem}
\begin{proof}
(a) It follows from Lemma \ref{L:weakly_conv} that $\Delta \overline{u}^{(n)}\to\Delta u_m$ on $\Omega_m^\circ\times[0,T]$.
Using \eqref{e:bounded_1}, we get $|\overline{u}|\le C_5/{\sqrt{\mu_0}}$, and hence $|\Delta \overline{u}^{(n)}|\le 2D_\mu\cdot C_5/{\sqrt{\mu_0}}$ on $\Omega_m^\circ\times [0,T]$. Combining this with Dominated Convergence Theorem, we get the result hold.

(b) From \eqref{e:inner_conve_unk_to_u}, we have
\begin{eqnarray*}
|\overline{u}^{(n)}|^{p-1}\cdot \overline{u}^{(n)}\to |u_m|^{p-1}\cdot u_m\quad\mbox{on }\Omega_m^\circ\times[0,T].
\end{eqnarray*}
Obviously, $\|\overline{u}^{(n)}(\cdot,t)\|_{L^{2p}(\Omega_m)}\le C_5$ implies $|\overline{u}|^{2p}\le C_5^{2p}/{\mu_0}$ on $\Omega_m^\circ\times[0,T]$. Consequently, the result hold follows from Dominated Convergence Theorem.
\end{proof}

\begin{lem}\label{L:unique_solution_on_Omega_m_eq}
There is a unique solution $u_m\in\mathcal{H}_m$ of \eqref{e:heat_eq_on_graph_m}.
\end{lem}
\begin{proof} This proof consists two parts.

{\bf Part 1. Existence}
In view of \eqref{e:delta2_ui_is_a_uni_sol}--\eqref{e:defi_delta_overline_u(n)}, we get
\begin{eqnarray*}
\int_0^T\big[\partial_t u^{(n)}(x,t)+|\overline{u}^{(n)}(x,t)|^{p-1}\cdot\overline{u}^{(n)}(x,t)
-\Delta\overline{u}^{(n)}(x,t)\big]\,dt=0\qquad\mbox{on }\Omega_m^\circ.
\end{eqnarray*}
Taking the limit as $n\to\infty$, using Lemmas \ref{L:weakly_conv} and \ref{L:int_0_T_conv}, we get
\begin{eqnarray*}\label{e:0_t_conver}
\int_0^T\big[\partial_tu_m(x,t)+|u_m(x,t)|^{p-1}\cdot u_m(x,t)-\Delta u_m(x,t)\,dt=0\quad\mbox{on }\Omega_m^\circ.
\end{eqnarray*}
Also, $u_m(x,0)=h_m(x)$ on $\Omega^\circ_m$, and $u_m=0$ on $V\setminus\Omega^\circ_m \times[0,T]$. 
Note that $T$ is arbitrary, we get $u_m$ is a solution of \eqref{e:heat_eq_on_graph_m}.

{\bf Part 2. Uniqueness} Suppose that $\hat{u}$ is another solution of \eqref{e:heat_eq_on_graph_m}. Let $w:=u_m-\hat{u}$. Then
\begin{equation*}
\left\{
\begin{aligned}
&\partial_tw+|u_m|^{p-1}\cdot u_m-|\hat{u}|^{p-1}\cdot\hat{u}=\Delta w,\qquad&&(x,t)\in\Omega_m^\circ\times[0,\infty),\\
&w(x,0)=0,\qquad&&x\in\Omega_m^\circ,\\
&w(x,t)=0,\qquad&&(x,t)\in V\setminus\Omega_m^\circ\times[0,\infty).
\end{aligned}
\right.
\end{equation*}

Consider the following energy functional
\begin{eqnarray*}
E_m(t)=\int_{\Omega_m} |w(x,t)|^2\,d\mu\qquad\mbox{ for }t\ge 0\mbox{ and }m\ge1.
\end{eqnarray*}
Using \eqref{e:green_formula_eq}, we obtain
\begin{eqnarray*}
\begin{aligned}
E'_m(t)=&2\int_{\Omega_m} w(x,t)\cdot \partial_tw(x,t)\,d\mu\\
=&2\int_{\Omega_m^\circ}w(x,t)\cdot
\big[\Delta w(x,t)-\big(|u_m(x,t)|^{p-1}\cdot u_m(x,t)-|\hat{u}(x,t)|^{p-1}\cdot \hat{u}(x,t)\big)\big]\,d\mu\\
=&-2\int_{\Omega_m^\circ}\big(u_m(x,t)-\hat{u}(x,t)\big)
\big(|u_m(x,t)|^{p-1}\cdot u_m(x,t)-|\hat{u}(x,t)|^{p-1}\cdot \hat{u}(x,t)\big)\,d\mu\\
&-2\int_{\Omega_m}|\nabla w(x,t)|^2\,d\mu.
\end{aligned}
\end{eqnarray*}
Note that for any $p\ge 1$,
\begin{eqnarray*}\label{e:ineq_u_and_hat_u}
\big(u_m-\hat{u}\big)\cdot\big(|u_m|^{p-1}\cdot u_m-|\hat{u}|^{p-1}\cdot\hat{u}\big)\ge0\quad\mbox{ on }\Omega_m^\circ\times[0,\infty).
\end{eqnarray*}
So $E_m'(t)\le 0$, which, together with $E_m(0)=0$, yields $E_m(t)=0$ on $[0,\infty)$.
Hence $u_m=\hat{u}$.
\end{proof}

\subsection{Convergence of $u_m$.}\label{SS:convergence_um_eq}
In this subsection, we prove $u_m$ converges to a solution $u$ of \eqref{e:heat_eq_on_graph}.

Let $Y$ be a Hilbert space. Also, let
\begin{eqnarray*}
L^2([0,T]; Y):=\left\{v(t): [0,T]\to Y\Big|\int_0^T\|v(t)\|_Y\,dt<\infty \mbox{ in Bochner sense}\right\}.
\end{eqnarray*}
Then $L^2([0,T]; Y)$ is a Hilbert space with the following inner product
\begin{eqnarray*}
(v_1, v_2)_{L^2([0,T]; Y)}=\int_0^T (v_1(t), v_2(t))_Y\,dt,
\end{eqnarray*}
where $\|\cdot\|_Y$ and $(\cdot,\cdot)_Y$ are norm and inner product in $Y$, respectively. For $v(t)\in L^2([0,T]; Y)$, the integral $\int_0^t v(s)\,ds=\mathcal{V}(t)$, in the Bochner sense, of $v$ is defined by
\begin{eqnarray*}
(\mathcal{V}(t), \phi)_Y=\int_0^t(v(s), \phi)_Y\,ds\quad\mbox{ for any }\phi\in Y.
\end{eqnarray*}
It follows that $v(t)$ can be identified a continuous function that is mapped from $[0,T]$ to $Y$ if $v(t)\in L^2([0,T]; Y)$. Moreover, $\mathcal{V}'(t)=v(t)$ for almost all $t\in [0,T]$ in the sense
$$\lim_{\widetilde{t}\to t}\|\frac{\mathcal{V}(\widetilde{t})-\mathcal{V}(t)}{\widetilde{t}-t}-v(t)\|_Y=0.$$

\begin{lem}\label{L:uni_bdd_u_m_eq}
For any $m\ge1$ and any $T>0$, the following results hold:
\begin{enumerate}
\item[(a)] $u_m$ is uniformly bounded in $L^2([0,T]; \mathcal{H})$;
\item[(b)] $\partial_t u_m$ are uniformly bounded in $L^2([0,T]; \mathcal{L}^2)$.
\end{enumerate}
\end{lem}
\begin{proof}
(a) Lemma \ref{L:unique_solution_on_Omega_m_eq} implies that for any $v\in\mathcal{H}_m$,
\begin{equation}\label{e:eq_u_m_v}
\int_{\Omega^\circ}\partial_tu_m\cdot v\,d\mu+\int_{\Omega^\circ}|u_m|^{p-1}\cdot u_m\cdot v\,d\mu
=\int_{\Omega^\circ}\Delta u_m\cdot v\,d\mu.
\end{equation}

Letting $v=u_m$ in \eqref{e:eq_u_m_v}, we get
\begin{equation}\label{e:ineq_u_m_eq}
\frac{d}{dt}\|u_m(\cdot,t)\|^2_{L^2(\Omega^\circ)}
+2\|u_m(\cdot,t)\|^{p+1}_{L^{p+1}(\Omega^\circ)}+2\|\nabla u_m(\cdot,t)\|^2_{L^2(\Omega)}=0.
\end{equation}
Obviously,
\begin{eqnarray*}
\frac{d}{dt}\|u_m(\cdot,t)\|^2_{L^2(\Omega^\circ)}+2\|\nabla u_m(\cdot,t)\|^2_{L^2(\Omega)}\le0.
\end{eqnarray*}
Integrating from $0$ to $T$ on the both sides of the above inequality, we obtain
\begin{eqnarray*}\label{e:bdd_int_0_T_nabla_um_var}
\int_0^T\|\nabla u_m(\cdot,t)\|^2_{L^2(\Omega)}\,dt
\le\frac{1}{2}\|u_m(\cdot,0)\|^2_{L^2(\Omega^\circ)}-\frac{1}{2}\|u_m(\cdot,T)\|^2_{L^2(\Omega^\circ)}
\le\|h\|^2_{L^2(\Omega^\circ)}.
\end{eqnarray*}

In view of \eqref{e:ineq_u_m_eq}, we have $\frac{d}{dt}\|u_m(\cdot,t)\|^2_{L^2(\Omega^\circ)}\le 0$. It follows that
for any $t\in[0,T]$,
\begin{eqnarray}\label{e:L2_u_m_is_bdd_eq}
\|u_m(\cdot,t)\|^2_{L^2(\Omega)}\le \|h_m\|^2_{L^2(\Omega^\circ)}\le \|h\|^2_{L^2(\Omega^\circ)}.
\end{eqnarray}
Hence
\begin{eqnarray*}\label{e:L2_0T_L2_u_m_is_bdd_eq}
\int_0^T\|u_m(\cdot,t)\|^2_{L^2(\Omega)}\,dt\le T \|h\|^2_{L^2(\Omega)}.
\end{eqnarray*}
Furthermore,
\begin{eqnarray*}\label{e:H_0T_L2_u_m_is_bdd_eq}
\int_0^T\|u_m(\cdot,t)\|^2_{\mathcal{H}}\,dt\le (T+1)\|h\|^2_{L^2(\Omega^\circ)}.
\end{eqnarray*}

(b) The fact $|\Delta u_m|^2\le D_\mu |\nabla u_m|^2$ on $\Omega_m^\circ\times[0,T]$ implies that
\begin{eqnarray*}
\int_0^T\|\Delta u_m(\cdot,t)\|^2_{L^2(\Omega^\circ)}\le D_\mu \|h\|^2_{L^2(\Omega^\circ)}.
\end{eqnarray*}

Letting $v=|u_m|^{p-1}\cdot u_m$ in \eqref{e:eq_u_m_v}, and using H\"older inequality, we see that
\begin{eqnarray*}
\frac{2}{p+1}\cdot\frac{d}{dt}\|u_m(\cdot,t)\|^{p+1}_{L^{p+1}(\Omega)}+\|u_m(\cdot,t)\|^{2p}_{L^{2p}(\Omega)}
\le\|\Delta u_m(\cdot,t)\|^2_{L^2(\Omega^\circ)}.
\end{eqnarray*}
Integrating from $0$ to $T$, we get
\begin{eqnarray*}\label{e:L2p_0T_L2p_u_m_is_bdd_eq}
\begin{aligned}
\int_0^T\|u_m(\cdot,t)\|^{2p}_{L^{2p}(\Omega)}\,dt
\le&\frac{2}{p+1}\|h\|^{p+1}_{L^{p+1}(\Omega)}+D_\mu\|h\|^2_{L^2(\Omega)}.
\end{aligned}
\end{eqnarray*}

Putting $v=\partial_tu_m$ in \eqref{e:eq_u_m_v}, using H\"older inequality and $(a-b)^2\le 2(a^2+b^2)$, we get
\begin{eqnarray*}
\int_{\Omega}|\partial_tu_m|^2\,d\mu\le 2\int_{\Omega^\circ}|\Delta u_m|^2\,d\mu+2\int_{\Omega}|u_m|^{2p}\,d\mu.
\end{eqnarray*}
Integrating from $0$ to $T$, we have
\begin{eqnarray*}
\int_0^T\|\partial_tu_m(\cdot,t)\|^2_{L^2(\Omega)}\,dt
\le\frac{4}{p+1}\|h\|^{p+1}_{L^{p+1}(\Omega^\circ)}+4D_\mu\|h\|^2_{L^2(\Omega^\circ)}.
\end{eqnarray*}
This completes the proof.
\end{proof}

In the rest of this subsection, we assume that $h\in L^2(\Omega^\circ)$ holds.
\begin{lem}\label{L:point_conve_eq}
There exist a function $u\in \mathcal{L}^2$ satisfying $\partial_t u\in \mathcal{L}^2$ and a subsequence of $\{u_{m}\}$, denoted by itself, such that the following results hold on $\Omega^\circ\times [0,T]$:
\begin{enumerate}
\item[(a)] $u_{m}\to u$ and $\Delta u_{m}\to \Delta u$;
\item[(b)] $\partial_t u_{m}\to \partial_t u$;
\item[(c)] $|u_{m}|^{p-1}\cdot u_{m}\to |u|^{p-1}\cdot u$.
\end{enumerate}
\end{lem}
\begin{proof}
(a) Lemma \ref{L:mathcal_L_and_Lm_are_closed} and \eqref{e:L2_u_m_is_bdd_eq} imply that there exists a function $u\in\mathcal{L}^2$ and a subsequence of $\{u_{m}\}$, denoted by itself, such that
\begin{eqnarray}\label{e:conve_inner_L2_eq}
(u_{m}(\cdot,t), v)\to (u(\cdot,t), v) \quad\mbox{ for any }t\in[0,T]\mbox{ and any }v\in\mathcal{L}^2.
\end{eqnarray}
For any $x_0\in \Omega^\circ$ and $t\in[0,T]$, define
\begin{equation*}\label{e:defi_delta_eq}
\delta_{x_0}(x,t):=
\left\{
\begin{aligned}
&1,\quad&&x=x_0,\\
&0,\quad&&x\neq x_0.
\end{aligned}
\right.
\end{equation*}
Letting $v=\delta_{x_0}$ in \eqref{e:conve_inner_L2_eq}, we obtain
\begin{equation*}\label{e:point_conver_eq}
u_{m}(x_0,t)\to u(x_0,t)\qquad\mbox{on }[0,T].
\end{equation*}
The arbitrary of $x_0$ implies that $u_{m}\to u$ and $\Delta u_{m}\to \Delta u$ on $\Omega^\circ\times[0,T]$.

(b) The proof of this part consists two steps.

{\bf Step 1.} We prove $\partial_t u$ exists. Lemmas \ref{L:uni_bdd_u_m_eq}(b) and \ref{L:mathcal_L_and_Lm_are_closed} imply that there exists a subsequence of $\{\partial_t u_{m}\}$, denoted by itself, and a function $S\in L^2([0,T]; \mathcal{L}^2)$ such that
\begin{eqnarray*}
\int_0^T(\partial_t u_{m}-S, v)\,dt\to 0\qquad\mbox{ for any }v\in L^2([0,T];\mathcal{L}^2).
\end{eqnarray*}

In view of
\begin{eqnarray*}
u_{m}(x,T)-h_{m}(x)=\int_0^T\partial_t u_{m}(x,t)\,dt\quad\mbox{on }\Omega,
\end{eqnarray*}
we have
\begin{eqnarray}\label{e:partial_s_u_m_ell_eq}
(u_{m}(\cdot,T)-h_{m},v)=\int_0^T\big( \partial_t u_{m}(\cdot,t),v\big)\,dt\qquad\mbox{for any }v\in \mathcal{L}^2.
\end{eqnarray}

For any $t\in[0,T]$, since $|h_{m}(x) v(x,t)|\le |h(x)||v(x,t)|$ and
\begin{eqnarray*}
\sum_{x\in\Omega^\circ}|h(x)||v(x,t)|\mu(x)
\le \Big(\sum_{x\in\Omega^\circ}|h(x)|^2\mu(x)\Big)^{1/2}\Big(\sum_{x\in\Omega^\circ}|v(x,t)|^2\mu(x)\Big)^{1/2}
<\infty,
\end{eqnarray*}
using the Weierstrass M-test, we see that $\sum_{x\in\Omega^\circ}h_{m}(x) v(x,t)\mu(x)$ converges uniformly with respect to $m$.
This leads to
\begin{eqnarray*}
\lim_{m\to\infty}(h_{m},v)
=\lim_{m\to\infty}\sum_{x\in\Omega^\circ} h_{m}(x) v(x,t) \mu(x)
=\sum_{x\in\Omega^\circ}h(x) v(x,t)\mu(x)=(h,v).
\end{eqnarray*}

Taking limit as $m\to\infty$ in \eqref{e:partial_s_u_m_ell_eq}, using \eqref{e:conve_inner_L2_eq}, we obtain
\begin{eqnarray*}
\begin{aligned}
(u(\cdot,T)-h,v)=\int_0^T(S(\cdot,t),v)\,dt=\Big(\int_0^TS(\cdot,t)\,dt,v\Big).
\end{aligned}
\end{eqnarray*}
In the second equality of the above equation, we use the fact that
\begin{equation}\label{e:fact_2_eq}
\int_0^T\Big(\sum_{x\in\Omega^\circ}S(x,t)v(x,t)\mu(x)\Big)\,dt=\sum_{x\in\Omega^\circ}\Big(\int_0^TS(x,t)\,dt\Big)v(x,t)\mu(x).
\end{equation}
The derivation of the above fact is as follows.
Since $S, v\in L^2([0,T],\mathcal{L}^2)$, that is,
$$\int_0^T\sum_{x\in\Omega^\circ}|S(x,t)|^2\mu(x)\,dt<\infty\quad\mbox{ and }\quad\int_0^T\sum_{x\in\Omega^\circ}|v(x,t)|^2\mu(x)\,dt<\infty,$$
we get $\sum_{x\in\Omega^\circ}\big|S(x,t)\sqrt{\mu(x)}\big|^2$ and $\sum_{x\in\Omega^\circ}\big|v(x,t)\sqrt{\mu(x)}\big|^2$ are convergent on $[0,T]$, and hence $\sum_{x\in\Omega^\circ}|S(x,t)v(x,t)\mu(x)|$ is convergent on $[0,T]$.
This leads to $\sum_{x\in\Omega^\circ}S(x,t)v(x,t)\mu(x)$ is convergent on $[0,T]$. Also, $\sum_{x\in\Omega^\circ}S(x,t)v(x,t)\mu(x)$ is integrated on $[0,T]$, since
\begin{eqnarray*}
\begin{aligned}
\int_0^T\sum_{x\in\Omega^\circ}S(x,t)v(x,t)\mu(x)\,dt\le \frac{1}{2}\Big(\int_0^T\sum_{x\in\Omega^\circ}S^2(x,t)\mu(x)\,dt+\int_0^T\sum_{x\in\Omega^\circ}v(x,t)^2\mu(x)\,dt\Big)
<\infty.
\end{aligned}
\end{eqnarray*}
Consequently,
\begin{eqnarray*}
\int_0^T\Big(\sum_{x\in\Omega^\circ}S(x,t)v(x,t)\mu(x)\Big)\,dt
=\sum_{x\in\Omega^\circ}\int_0^T\Big(S(x,t)v(x,t)\mu(x)\Big)\,dt.
\end{eqnarray*}

By the above derivation, we get $S$ is integrate on $[0,T]$, $u$ is continuous with respect to $t$ and $u(x,t)-h(x)=\int_0^tS(x,s)\,ds$ for any $t\in[0,T]$. Moreover, $u(x,0)=h(x)$ and $S=\partial_t u$.

{\bf Step 2.} We show that $\partial_t u_{m}\to \partial_t u$. For any $t_0\in[0,T]$ and $x\in\Omega$,
\begin{eqnarray*}
\begin{aligned}
\lim_{m\to\infty}\partial_t u_{m}(x,t_0)
=&\lim_{m\to\infty}\lim_{t\to t_0}\frac{u_{m}(x,t)-u_{m}(x,t_0)}{t-t_0}\\
=&\lim_{t\to t_0}\lim_{m\to\infty}\frac{u_{m}(x,t)-u_{m}(x,t_0)}{t-t_0}\\
=&\lim_{t\to t_0}\frac{u(x,t)-u(x,t_0)}{t-t_0}\\
=&\partial_t u(x,t_0).
\end{aligned}
\end{eqnarray*}
This proves (b).

(c) Since $u_{m}$ and $u$ are bounded on $\Omega^\circ\times[0,T]$, using (a), we get
\begin{eqnarray*}
\begin{aligned}
|u_{m}|^{p-1}-|u|^{p-1}
=\big(|u_{m}|-|u|\big)\cdot\big(|u_{m}|^{p-2}+|u_{m}|^{p-3}\cdot|u|+\cdots+|u|^{p-2}\big)
\to0,
\end{aligned}
\end{eqnarray*}
that is $|u_{m}|^{p-1}\to |u|^{p-1}$ for $p\ge1$. By triangle inequality, we obtain
\begin{eqnarray*}
\big||u_{m}|^{p-1}\cdot u_{m}-|u|^{p-1}\cdot u\big|\to 0.
\end{eqnarray*}
This completes the proof.
\end{proof}

Under the above derivation, we give the proof of Theorem~\ref{T:sol_of_heat_equ}.
\begin{proof}[Proof of Theorem~\ref{T:sol_of_heat_equ}]
Lemma \ref{L:unique_solution_on_Omega_m_eq} implies that for any $T>0$,
\begin{equation*}
\left\{
\begin{aligned}
&\partial_t u_m+|u_m|^{p-1}\cdot u_m=\Delta u_m,\quad&&(x,t)\in\Omega_m^\circ\times[0,T],\\
&u_m(x,0)=h_m(x),\quad &&x\in\Omega_m^\circ,\\
&u_m(x,t)=0,\quad&&(x,t)\in V\setminus\Omega_m^\circ\times[0,T].
\end{aligned}
\right.
\end{equation*}
Using Lemma \ref{L:point_conve_eq}, we see that there exists a function $u(\cdot,t)\in\mathcal{L}^2$ such that
$$\partial_t u+u\cdot|u|^{p-1}=\Delta u\qquad\mbox{on }\Omega^\circ\times[0,T],$$
$u(x,0)=h(x)$ on $\Omega^\circ$, and $u(x,t)=0$ on $V\setminus\Omega^\circ\times[0,T]$. 
The result follows $T$ is arbitrary.
\end{proof}

\section{Proof of Corollary~\ref{C:sol_of_heat_equ_p1}}\label{S:proof_thm_heat_equ_p1}
\setcounter{equation}{0}

\begin{proof}[Proof of Corollary~\ref{C:sol_of_heat_equ_p1}]
Let $u(x,t)=\sum\limits_{j=1}^N h_j(t)\varphi_j(x)$ for $N=\#\Omega^\circ$. Then
$$h(x)=\sum_{j=1}^Nh_j(0)\varphi_j(x),\quad \partial_tu(x,t)=\sum_{j=1}^N \dot{h}_j(t)\varphi_j(x),\quad
\mbox{and }-\Delta u(x,t)=\sum_{j=1}^N\lambda_j h_j(t)\varphi_j(x),$$
where $\cdot$ denote the derivative with respect to $t$.
For $p=1$, the first equation in \eqref{e:heat_eq_on_finite_graph} becomes
\begin{eqnarray*}
\sum_{j=1}^N \dot{h}_j(t)\varphi_j(x)+\sum_{j=1}^N h_j(t)\varphi_j(x)+\sum_{j=1}^N \lambda_j h_j(t)\varphi_j(x)=0.
\end{eqnarray*}
This leads to for any $j=1,\ldots, N$,
\begin{eqnarray*}
\dot{h}_j(t)+h_j(t)+\lambda_j h_j(t)=0.
\end{eqnarray*}
Hence $h_j(t)=h_j(0)e^{-(\lambda_j+1)t}$, and so
\begin{eqnarray*}
u(x,t)=\sum\limits_{j=1}^N h_j(0)e^{-(\lambda_j+1)t} \varphi_j(x),
\end{eqnarray*}
where $h_j(0)=(h, \varphi_j)_{W_0^{1,2}(\Omega)}$. This completes the proof.
\end{proof}

\section{Parabolic variational inequalities}\label{S:para_var_ineq}
\setcounter{equation}{0}

In this section, we study the existence of a solution of \eqref{e:para_var_ineq}--\eqref{e:initial_boundary_value_var}.
Let $G$ be a locally finite connected weighted graph, and $\Omega\subset V$ be unbounded. Suppose there exists a sequence of bounded subsets $\{\Omega_m\}_{m\ge1}$ of $\Omega$ satisfying \eqref{e:construction_of_Omega_m_eq}. Let $\mathcal{L}^2_m,\ \mathcal{L}^2,\ \mathcal{H}_m$ and $\mathcal{H}$ be defined as in subsection \ref{SS:constrction_of_bdd_subsets_eq}. For given functions $f$ and $g$,
let $f_m:=f|_{\Omega_m^\circ\times[0,\infty)}$ and $g_m:=g|_{\Omega_m^\circ}$. As same as Section \ref{S:proof_thm_heat_equ}, we divide the research into two steps (see subsection \ref{SS:par_var_ine_on_Omega_m_var}--\ref{SS:convergence_um_var}).

Unless otherwise specified, the symbols in this Section have nothing to do with these in Section \ref{S:proof_thm_heat_equ}.

\subsection{Parabolic variational inequality on $\Omega_m$}\label{SS:par_var_ine_on_Omega_m_var}

In this subsection, using Rothe's method, we study the following parabolic variational inequality
\begin{eqnarray}\label{e:para_var_ineq_m}
\int_{\Omega_m^\circ} \partial_t u\cdot(v-u)\,d\mu\ge \int_{\Omega_m^\circ}(\Delta u+f_m)\cdot(v-u)\,d\mu
\quad\mbox{ for }t\in [0,\infty)\mbox{ and any }v\in \mathcal{H}_m,
\end{eqnarray}
with
\begin{eqnarray}\label{e:initial_boundary_value_var_m}
u(x,0)=g_m(x)\quad\mbox{ for }x\in \Omega_m^\circ,
\qquad\mbox{and}\qquad u(x,t)=0\quad\mbox{ for }(x,t)\in V\setminus\Omega_m^\circ\times[0,\infty).
\end{eqnarray}

For any $T>0$, we divide $[0,T]$ into $n$ subintervals $[t_{i-1}, t_i]$ satisfying $t_i-t_{i-1}=\ell$ for $i\in \Lambda:=\{1,\ldots, n\}$, $t_0=0$, and $n\ell=T$. Define $u^m_{n,0}(x):=g_m(x)$, and
\begin{equation*}\label{e:defi_fi_var}
f^m_{n,i}(x):=f_m(x,t_i)\qquad\mbox{for }i\in \{0\}\cup \Lambda.
\end{equation*}
Successively, for $i\in \Lambda$, we want to find a unique $u^{m}_{n,i}\in \mathcal{H}_m$ satisfies \begin{equation}\label{e:basic_ineq_var}
\int_{\Omega_m^\circ}(u-u_{n,i-1})/\ell\cdot(v-u)\,d\mu\ge \int_{\Omega_m^\circ}(\Delta u+f^m_{n,i})\cdot(v-u)\,d\mu
\quad\mbox{ for any }v\in \mathcal{H}_m.
\end{equation}
For $m\ge1$ and $i\in\{0\}\cup\Lambda$, we will denote $u^m_{n,i}$ and $f^m_{n,i}$ simply by $u_{n,i}$ and $f_{n,i}$, respectively.

\begin{lem}\label{L:uni_sol_var}
For $m\ge1$ and any $i\in \Lambda$, \eqref{e:basic_ineq_var} has a unique solution $u_{n,i}\in \mathcal{H}_m$.
\end{lem}

\begin{proof}
If $i=1$, then \eqref{e:basic_ineq_var} becomes
\begin{equation}\label{e:proof_Lemma_uni_sol_var_1}
\int_{\Omega_m^\circ}(u/\ell-\Delta u)\cdot (v-u)\,d\mu\ge \int_{\Omega_m^\circ}(f_{n,1}+g_{m}/\ell)\cdot(v-u)\,d\mu
\quad\mbox{ for any }v\in \mathcal{H}_m.
\end{equation}
It is obvious that $\mathcal{H}_m$ is a closed and convex set.
Next, We show that \eqref{e:proof_Lemma_uni_sol_var_1} has a unique solution $u_{n,1}\in \mathcal{H}_m$.
Let $$a(v_1,v_2):=1/\ell\int_{\Omega_m^\circ} v_1\cdot v_2\,d\mu-\int_{\Omega_m^\circ}\Delta v_1\cdot v_2\,d\mu\quad\mbox{ for }v_1,v_2\in \mathcal{H}_m.$$
Then $a(\cdot, \cdot)$ is a coercive bilinear form on $\mathcal{H}_m$.
In fact, using \eqref{e:green_formula_eq}, we get
\begin{eqnarray*}
\begin{aligned}
a(v,v)
=1/\ell\int_{\Omega_m^\circ} |v|^2\,d\mu+\int_{\Omega_m}|\nabla v|^2\,d\mu
\ge\beta\|v\|_{\mathcal{H}_m}^2\quad\mbox{ for any }v\in \mathcal{H}_m,
\end{aligned}
\end{eqnarray*}
where $\beta:=\min\{1/\ell,1\}$. Let $F(u):=\int_{\Omega_m^\circ}(f_{n,1}+g_m/\ell)\cdot u\,d\mu$ for $u\in \mathcal{H}_m$.
Since $f_m$ and $g_m$ are given, we get $F$ is linear and bounded on $\mathcal{H}_m$.
It follows from Theorem \ref{T:Kinderlehrer-Stampacchia} that \eqref{e:proof_Lemma_uni_sol_var_1} has a unique solution $u_{n,1}$ in $\mathcal{H}_m$.

Similarly, there exists a unique solution $u_{n,i}\in \mathcal{H}_m$ of \eqref{e:basic_ineq_var} for $i=\Lambda\setminus\{1\}$.
\end{proof}

Let $u_{n,i}(x)$ be the approximation of a solution $u_m(x,t)$ for \eqref{e:para_var_ineq_m}--\eqref{e:initial_boundary_value_var_m} at $t=t_i$. Also, let
\begin{equation}\label{e:delta_ui_var}
\delta u^m_{n,i}(x)=(u_{n,i}(x)-u_{n,i-1}(x))/\ell\qquad \mbox{ for }i\in \Lambda,
\end{equation}
and denote simply $\delta u^m_{n,i}$ by $\delta u_{n,i}$. Define Rothe's function $u^{(n),m}(x,t)$ (or simply $u^{(n)}(x,t)$) from $[0,T]$ to $\mathcal{H}_m$ by
\begin{equation}\label{e:defi_un_and_fn_var}
u^{(n)}(x,t)=u_{n,i-1}(x)+(t-t_{i-1})\cdot\delta u_{n,i}(x).
\end{equation}
We also define step functions
\begin{equation}\label{e:defi_overline_u(n)_var}
\overline{u}^{(n)}(x,t)=
\left\{
\begin{aligned}
&u_{n,i}(x),\qquad &&(x,t)\in\Omega_m^\circ\times(t_{i-1},t_i],\\
&g_m(x),\qquad &&(x,t)\in\Omega_m^\circ\times[-\ell,0],\\
&0,\qquad&&(x,t)\in V\setminus\Omega_m\times[-\ell,0],
\end{aligned}
\right.
\end{equation}
and
\begin{equation}\label{e:defi_overline_f(n)_var}
\overline{f}^{(n)}(x,t)=
\left\{
\begin{aligned}
&f_{n,i}(x),\qquad&& (x,t)\in\Omega_m^\circ\times(t_{i-1},t_i],\\
&f_{n,0}(x),\qquad&& (x,t)\in\Omega_m^\circ\times[-\ell,0].
\end{aligned}
\right.
\end{equation}
Based on the above derivation, we can get the following estimates on $u_{n,i}$ and $\delta u_{n,i}$. From now on, we assume that \eqref{e:condition_for_f} holds.
\begin{lem}\label{L:est_delta_ui_var}
There exist positive constants $C_1,\ C_2$, depending only on $\Omega_m$ and $T$,
such that
\begin{eqnarray*}
\|\delta u_{n,i}\|_{L^2(\Omega_m)}\le C_1,\qquad \|u_{n,i}\|_{L^2(\Omega_m)}\le C_2\quad\mbox{ for any }i\in \Lambda.
\end{eqnarray*}
\end{lem}

\begin{proof}
It follows from \eqref{e:basic_ineq_var} and Lemma \ref{L:uni_sol_var} that for any $i\in \Lambda$ and any $v\in \mathcal{H}_m$,
\begin{equation}\label{e:basic_ineq_var_2}
\int_{\Omega_m^\circ}(u_{n,i}-u_{n,i-1})/\ell\cdot(v-u_{n,i})\,d\mu
\ge \int_{\Omega_m^\circ}(\Delta u_{n,i}+f_{n,i})\cdot(v-u_{n,i})\,d\mu.
\end{equation}

Letting $i=j$, $v=u_{n,j-1}$ and $i=j-1$, $v=u_{n,j}$ in \eqref{e:basic_ineq_var_2}, and adding those inequalities, we obtain
\begin{eqnarray*}
\begin{aligned}
&\int_{\Omega_m^\circ}(u_{n,j}-2u_{n,j-1}+u_{n,j-2})/\ell\cdot(u_{n,j}-u_{n,j-1})\,d\mu\\
\le&\int_{\Omega_m^\circ}\big(\Delta(u_{n,j}-u_{n,j-1})+(f_{n,j}-f_{n,j-1})\big)\cdot(u_{n,j}-u_{n,j-1})\,d\mu.
\end{aligned}
\end{eqnarray*}
Multiplying both sides of the above inequality by $1/\ell$, using \eqref{e:green_formula_eq}, we get
\begin{eqnarray*}
\begin{aligned}
&\int_{\Omega_m^\circ}\big|(u_{n,j}-u_{n,j-1})/\ell\big|^2\,d\mu\\
\le&\int_{\Omega_m^\circ}\big|(u_{n,j}-u_{n,j-1})/\ell\big|^2\,d\mu-1/\ell\cdot\int_{\Omega_m^\circ}\Delta(u_{n,j}-u_{n,j-1})\cdot(u_{n,j}-u_{n,j-1})\,d\mu\\
\le&\Big[\Big(\int_{\Omega_m^\circ}\big|(u_{n,j-1}-u_{n,j-2})/\ell\big|^2\,d\mu\Big)^{1/2}
+\Big(\int_{\Omega_m^\circ}|(f_{n,j}-f_{n,j-1})|^2\,d\mu\Big)^{1/2}\Big]\\
&\cdot\Big(\int_{\Omega_m^\circ}\big|(u_{n,j}-u_{n,j-1})/\ell\big|^2\,d\mu\Big)^{1/2},
\end{aligned}
\end{eqnarray*}
which yields
\begin{eqnarray}\label{e:ineq_delta_ui_recurrence}
\|(u_{n,j}-u_{n,j-1})/\ell\|_{L^2(\Omega_m^\circ)}
\le \|(u_{n,j-1}-u_{n,j-2})/\ell\|_{L^2(\Omega_m^\circ)}+\|f_{n,j}-f_{n,j-1}\|_{L^2(\Omega_m^\circ)}.
\end{eqnarray}

Similarly, letting $i=1$ and $v=u_{n,0}$ in \eqref{e:basic_ineq_var_2}, we get
$$\big\|(u_{n,1}-u_{n,0})/\ell\big\|_{L^2(\Omega_m^\circ)}
\le\|\Delta u_{n,0}+f_{n,1}\|_{L^2(\Omega_m^\circ)}
\le\|\Delta g_m\|_{L^2(\Omega_m)}+\|f_{n,1}\|_{L^2(\Omega_m^\circ)},$$
which, together with \eqref{e:ineq_delta_ui_recurrence} and \eqref{e:condition_for_f}, yields
\begin{eqnarray*}
\begin{aligned}
\|(u_{n,i}-u_{n,i-1})/\ell\|_{L^2(\Omega_m^\circ)}
\le&\|\Delta g_m\|_{L^2(\Omega_m^\circ)}+\|f_{n,0}\|_{L^2(\Omega_m^\circ)}
+\sum_{k=0}^{i-1}\big\|f_{n,k+1}-f_{n,k}\big\|_{L^2(\Omega_m^\circ)}\\
\le&\|\Delta g_m\|_{L^2(\Omega_m^\circ)}+\|f(\cdot,0)\|_{L^2(\Omega_m^\circ)}+c\cdot n\cdot T/n\\
\le& C_1\qquad\mbox{for some constant } C_1:=C_1(\Omega_m, T)>0.
\end{aligned}
\end{eqnarray*}
Hence for any $i\in \Lambda$, $\|\delta u_{n,i}\|_{L^2(\Omega_m^\circ)}\le C_1$, that is, $\|u_{n,i}-u_{n,i-1}\|_{L^2(\Omega_m^\circ)}\le C_1\ell$. Thus for any $i\in \Lambda$,
\begin{eqnarray*}
\begin{aligned}
\|u_{n,i}\|_{L^2(\Omega_m^\circ)}
\le\|u_{n,0}\|_{L^2(\Omega_m^\circ)}+\sum_{k=0}^{i-1}\|u_{n,k+1}-u_{n,k}\|_{L^2(\Omega_m^\circ)}
\le\|g_m\|_{L^2(\Omega_m^\circ)}+C_1 T\le C_2,
\end{aligned}
\end{eqnarray*}
where $C_2:=C_2(\Omega_m, T)>0$ is a constant. This completes the proof.
\end{proof}

From Lemma \ref{L:est_delta_ui_var}, we can derived directly the following priori estimates.
\begin{lem}\label{L:est_of_u(n)_var}
There exist two positive constants $C_i:=C_i(\Omega_m, T)$, $i=3,4$, such that for any $t\in [0,T]$,
\begin{eqnarray*}
\|u^{(n)}(\cdot,t)\|_{L^{2}(\Omega_m)}+\|\overline{u}^{(n)}(\cdot,t)\|_{L^{2}(\Omega_m)}
+\|\partial_t u^{(n)}(\cdot,t)\|_{L^2(\Omega_m)}\le C_3,
\end{eqnarray*}
and
\begin{eqnarray*}
\|u^{(n)}(\cdot,t)-\overline{u}^{(n)}(\cdot,t)\|_{L^2(\Omega_m)}\le C_4/n.
\end{eqnarray*}
\end{lem}

Lemma \ref{L:est_of_u(n)_var} can deduce the following results. The proof is the same as that of Lemma \ref{L:weakly_conv}.
\begin{lem}\label{L:weakly_conv_var}
There exist a function $u_m\in \mathcal{L}^2_m$ satisfying $\partial_t u_m\in \mathcal{L}^2_m$,
and two subsequences of $\{u^{(n)}\}$, $\{\overline{u}^{(n)}\}$, denoted by themselves, such that on $\Omega^\circ_m\times[0,T]$,
\begin{enumerate}
\item[(a)] $u^{(n)}\to u_m$ and $\overline{u}^{(n)}\to u_m$;
\item[(b)] $\partial_t u^{(n)}\to \partial_t u_m$.
\end{enumerate}
\end{lem}

\begin{proof}
(a) follows from Lemma \ref{L:est_of_u(n)_var}.

(b) Using Lemma \ref{L:est_of_u(n)_var} again, we can extract a subsequence $\{\partial_t u^{(n)}\}$ such that
\begin{eqnarray*}
\partial_t u^{(n)}\to \partial_t u_m \quad\mbox{ on }\Omega^\circ_m\times[0,T],
\end{eqnarray*}
with $u_m(x,0)=g_m$. In the proof, we use the fact that
\begin{eqnarray*}
\int_0^t \partial_s u^{(n)}(x,s)\,ds\to\int_0^t \partial_t u_m(x,s)\,ds\quad\mbox{ on }\Omega^\circ_m\times[0,T].
\end{eqnarray*}
The proof is completed.
\end{proof}

\begin{lem}\label{L:int_int_conv_var}
The following results hold:
\begin{enumerate}
\item[(a)] $\partial_t u^{(n)} \to \partial_t u_m$ and $\Delta \overline{u}^{(n)}\to \Delta u_m$ in $L^2([0,T];\mathcal{L}^2_m)$;
\item[(b)] $\overline{f}^{(n)}\to f_m$ in $L^2([0,T]; L^{2}(\Omega_m^\circ))$;
\item[(c)] $\int_0^T\int_{\Omega_m^\circ} \partial_t u^{(n)}\cdot \overline{u}^{(n)} \,d\mu\,dt\to \int_0^T\int_{\Omega_m^\circ} \partial_t u_m\cdot u_m\,d\mu\,dt$;
\item[(d)] $\int_0^T\int_{\Omega_m}|\nabla\overline{u}^{(n)}|^2\,d\mu\,dt\to \int_0^T\int_{\Omega_m}|\nabla u_m|^2 \,d\mu\,dt$;
\item[(e)] $\int_0^T\int_{\Omega_m^\circ} \overline{f}^{(n)}\cdot \overline{u}^{(n)}\,d\mu\,dt\to
           \int_0^T\int_{\Omega_m^\circ}f_m\cdot u_m\,d\mu\,dt$.
\end{enumerate}
\end{lem}

\begin{proof}
(a) It follows from Lemma \ref{L:weakly_conv_var} that for any $t\in[0,T]$, $\partial_t u^{(n)}(\cdot,t) \to \partial_t u_m(\cdot,t)$ in $\mathcal{L}^2_m$. Note that $\|\partial_t u^{(n)}(\cdot,t)\|_{L^{2}(\Omega_m)}\le C_3$ on $[0,T]$, using Fatou's Lemma, we get
\begin{eqnarray*}
\int_{\Omega_m}|\partial u^{(n)}(x,t)-\partial_t u_m(x,t)|^2\,d\mu
\le 2 \int_{\Omega_m}|\partial u^{(n)}(x,t)|^2\,d\mu+2\liminf_{n\to\infty}\int_{\Omega_m}|\partial_t u^{(n)}(x,t)|^2\,d\mu\le 4C_3.
\end{eqnarray*}
The result $\partial_t u^{(n)} \to \partial_t u_m$ in $L^2([0,T];\mathcal{L}^2_m)$ follows from Dominated Convergence Theorem.

The fact $\overline{u}^{(n)}\to u_m$ on $(x,t)\in\Omega_m^\circ\times[0,T]$ implies that $\Delta \overline{u}^{(n)}\to\Delta u_m$.
Also, $\|\Delta \overline{u}^{(n)}\|_{L^{2}(\Omega^\circ_m)}\le 2D_\mu C_3$. Similarly, we get $\Delta \overline{u}^{(n)}\to \Delta u_m$ in $L^2([0,T];\mathcal{L}^2_m)$. This proves (a).

(b) Since $f_m$ is given on $\Omega_m^\circ$, we have $\|\overline{f}^{(n)}(\cdot, t)\|_{L^{2}(\Omega_m^\circ)}$ is uniformly bounded on $[0,T]$. It is easy to see that $\overline{f}^{(n)}\to f_m$ in $L^2([0,T]; L^{2}(\Omega_m^\circ))$.

(c) By Lemma \ref{L:weakly_conv_var} and triangle inequality, we get 
$\partial_t u^{(n)}\cdot \overline{u}^{(n)}\to \partial_t u_m\cdot u_m$ on $\Omega_m^\circ\times[0,T]$. In view of
\begin{eqnarray}\label{e:bounded_int_over_u_var}
\|\overline{u}^{(n)}(\cdot,t)\|_{L^{2}(\Omega_m)}+\|\partial_t u^{(n)}(\cdot,t)\|_{L^2(\Omega_m)}\le C_3,
\end{eqnarray}
we obtain
\begin{eqnarray*}
|\partial_t u^{(n)}|\le C_3/\sqrt{\mu_0}\quad\mbox{ and }\quad|\overline{u}^{(n)}|\le C_3/\sqrt{\mu_0}
\quad\mbox{ on }\Omega_m^\circ\times[0,T].
\end{eqnarray*}
Using Dominated Convergence Theorem, we have
$$\int_{\Omega_m^\circ} \partial_t u^{(n)}\cdot \overline{u}^{(n)} \,d\mu\to \int_{\Omega_m^\circ} \partial_t u_m\cdot u_m\,d\mu.$$
Using \eqref{e:bounded_int_over_u_var} and Dominated Convergence Theorem again, we get the result.

(d) In view of Lemma \ref{L:weakly_conv_var}(a), we have $|\nabla \overline{u}^{(n)}|^2\to |\nabla u_m|^2$ on $\Omega_m\times[0,T]$. Since
\begin{eqnarray*}
\begin{aligned}
\|\nabla \overline{u}^{(n)}(\cdot,t)\|^2_{L^2(\Omega_m)}
&=\sum_{x\in\Omega_m}\frac{1}{2\mu(x)}\sum_{y\sim x}\omega_{xy}|\overline{u}^{(n)}(y,t)-\overline{u}^{(n)}(x,t)|^2 \,\mu(x)\\
&=\frac{1}{2}\sum_{x\in\Omega_m}\sum_{y\sim x}\omega_{xy}|\overline{u}^{(n)}(y,t)-\overline{u}^{(n)}(x,t)|^2\\
&\le \sum_{x\in\Omega_m}\sum_{y\sim x}\omega_{xy}\big(|\overline{u}^{(n)}(y,t)|^2+|\overline{u}^{(n)}(x,t)|^2\big)\\
&\le M_d\sum_{x\in\Omega_m}\frac{1}{\mu(x)}\sum_{y\sim x}\omega_{xy}|\overline{u}^{(n)}(x,t)|^2\mu(x)\\
&\le D_\mu M_d^2\|\overline{u}^{(n)}(\cdot,t)\|^2_{L^2(\Omega_m)},
\end{aligned}
\end{eqnarray*}
using Dominated Convergence Theorem, we get the result holds.

(e) The proof is the same as that of (c).
\end{proof}

\begin{lem}\label{L:unique_solution_on_Omega_m}
There is a unique solution $u_m\in\mathcal{H}_m$ for \eqref{e:para_var_ineq_m} and \eqref{e:initial_boundary_value_var_m}.
\end{lem}

\begin{proof}
It follows from \eqref{e:basic_ineq_var}, \eqref{e:delta_ui_var}--\eqref{e:defi_overline_f(n)_var} and Lemma \ref{L:uni_sol_var}
that
\begin{eqnarray*}
\int_0^T\int_{\Omega_m^\circ}\big(\partial_t u^{(n)}(\cdot,t)-\Delta \overline{u}^{(n)}(\cdot,t)-\overline{f}^{(n)}(\cdot,t)\big)
\cdot\big(v-\overline{u}^{(n)}(\cdot,t)\big)\,d\mu\,dt\ge 0.
\end{eqnarray*}
Taking the limit as $n\to\infty$, using Lemma \ref{L:int_int_conv_var}, we obtain
\begin{eqnarray*}
\int_0^T\int_{\Omega_m^\circ}\big(\partial_t u_m-\Delta u_m-f_m\big)\cdot\big(v-u_m\big)\,d\mu\,dt\ge 0\qquad
\mbox{for any } v\in\mathcal{H}_m.
\end{eqnarray*}
In view of the arbitrary of $T$, we get
\begin{eqnarray*}
\int_{\Omega_m^\circ}\big(\partial_t u_m-\Delta u_m-f_m\big)\cdot\big(v-u_m\big)\,d\mu\ge 0\qquad \mbox{for any }v\in \mathcal{H}_m.
\end{eqnarray*}
From the proof of Lemma \ref{L:weakly_conv_var}, we get $u_m(x,0)=g_m(x)$ and $u_m\in \mathcal{H}_m$. Hence $u_m$ is a solution of \eqref{e:para_var_ineq_m}--\eqref{e:initial_boundary_value_var_m}.

Let $\widetilde{u}$ and $\hat{u}$ be two solutions of \eqref{e:para_var_ineq_m}--\eqref{e:initial_boundary_value_var_m}.
In \eqref{e:para_var_ineq_m}, letting $u=\widetilde{u}$, $v=\hat{u}$ and $u=\hat{u}$, $v=\widetilde{u}$, and
adding these inequalities, we obtain
\begin{eqnarray*}
\int_{\Omega_m^\circ}\partial_t(\widetilde{u}-\hat{u})\cdot(\hat{u}-\widetilde{u})\,d\mu
\ge \int_{\Omega_m^\circ}\Delta(\widetilde{u}-\hat{u})\cdot (\hat{u}-\widetilde{u})\,d\mu,
\end{eqnarray*}
and so
\begin{eqnarray*}
\begin{aligned}
\frac{1}{2}\frac{d}{dt}\|\hat{u}-\widetilde{u}\|^2_{L^2(\Omega_m^\circ)}
\le \int_{\Omega_m^\circ}\Delta(\hat{u}-\widetilde{u})\cdot (\hat{u}-\widetilde{u})\,d\mu
=-\int_{\Omega_m^\circ}|\nabla(\hat{u}-\widetilde{u})|^2\,d\mu
\le 0.
\end{aligned}
\end{eqnarray*}
On the other hand, $\frac{d}{dt}\|\hat{u}(\cdot,t)-\widetilde{u}(\cdot,t)\|^2_{L^2(\Omega_m^\circ)}\ge0$, which follows from
\begin{eqnarray*}
\int_0^t\frac{d}{ds}\|\hat{u}(\cdot,s)-\widetilde{u}(\cdot,s)\|^2_{L^2(\Omega_m^\circ)}\,ds
=\|\hat{u}(\cdot,t)-\widetilde{u}(\cdot,t)\|^2_{L^2(\Omega_m^\circ)}\ge0.
\end{eqnarray*}
Hence $\|\hat{u}(\cdot,t)-\widetilde{u}(\cdot,t)\|^2_{L^2(\Omega_m^\circ)}=0$. This leads to $\widetilde{u}=\hat{u}$.
\end{proof}

\subsection{Convergence of $u_m$}\label{SS:convergence_um_var}
In this subsection, we prove that $u_m$ converges to a solution $u$ of \eqref{e:para_var_ineq}--\eqref{e:initial_boundary_value_var}. In the rest of this subsection, we assume that $g\in L^2(\Omega^\circ)$.

\begin{lem}\label{L:uniform_bdd_seq_var}
$u_m$ and $\partial_t u_m$ are uniformly bounded in $\mathcal{L}^2$ and $L^2([0,T]; \mathcal{L}^2)$.
\end{lem}
\begin{proof}
{\bf Part 1.} $u_m$ is uniformly bounded in $\mathcal{L}^2$ and $L^2([0,T]; \mathcal{L}^2)$.

Lemma \ref{L:unique_solution_on_Omega_m} implies that for any $T>0$, any $t\in[0,T]$ and any $v\in\mathcal{H}_m$,
\begin{eqnarray}\label{e:um_on_Omega}
\int_{\Omega^\circ} \partial_t u_m\cdot (v-u_m)\,d\mu\ge \int_{\Omega^\circ} (\Delta u_m+f)\cdot (v-u_m)\,d\mu.
\end{eqnarray}

Letting $v=0$ in \eqref{e:um_on_Omega}, and using \eqref{e:green_formula_eq}, we obtain
\begin{eqnarray*}\label{e:ineq_dt_um_and_nabla_um_Omega}
\frac{d}{dt}\|u_m(\cdot,t)\|^2_{L^2(\Omega^\circ)}+2\|\nabla u_m(\cdot,t)\|^2_{L^2(\Omega)}
\le \|f(\cdot,t)\|^2_{L^2(\Omega^\circ)}+\|u_m(\cdot,t)\|^2_{L^2(\Omega^\circ)}.
\end{eqnarray*}
Let $\eta(t):=\|u_m(\cdot,t)\|^2_{L^2(\Omega^\circ)}$ and $\xi(t):=\|f(\cdot,t)\|^2_{L^2(\Omega^\circ)}$.
Then $\eta'(t)\le \eta(t)+\xi(t)$. This leads to
\begin{eqnarray*}
\eta(t)\le e^t\big(\eta(0)+\int_0^t\xi(s)e^{-s}\,ds\big).
\end{eqnarray*}
Hence
\begin{eqnarray*}
\begin{aligned}
\|u_m(\cdot,t)\|^2_{L^2(\Omega)}
\le &e^T\big(\|g\|^2_{L^2(\Omega^\circ)}+\|f(\cdot,t)\|^2_{L^2([0,T];\mathcal{L}^2)}\big)
\le  C\qquad\mbox{for any }t\in[0,T],
\end{aligned}
\end{eqnarray*}
where $C$ is a constant independent on $m$. Moreover, $\|u_m(\cdot,t)\|^2_{L^2([0,T];\mathcal{L}^2)}\le CT$. This proves that $u_m$ is uniformly bounded in $\mathcal{L}^2$ and $L^2([0,T];\mathcal{L}^2)$.

{\bf Part 2.} $\partial_t u_m$ is uniformly bounded in $\mathcal{L}^2$ and $L^2([0,T]; \mathcal{L}^2)$.

Next, we show that $u_m-\partial_t u_m\in\mathcal{H}_m$. Lemma \ref{L:weakly_conv_var} implies that $\partial_t u_m(\cdot, t)\in \mathcal{L}^2_m$ on $[0,T]$. Hence $\partial_t u_m$ is bounded on $\Omega_m^\circ\times[0,T]$, which yields
\begin{eqnarray*}
\|\nabla (\partial_t u_m)(\cdot,t)\|^2_{L^2(\Omega)}
=\|\nabla (\partial_t u_m)(\cdot,t)\|^2_{L^2(\Omega_m)}
\le(2D_\mu|\partial_t u_m(x,t)|)^2\cdot\mu(\Omega_m)<\infty,
\end{eqnarray*}
and so $\partial_t u_m(\cdot,t)\in W^{1,2}(\Omega_m)$ on $[0,T]$.
Since $u_m\in \mathcal{H}_m$, we get $u_m\equiv0$ on $V\setminus\Omega_m^\circ\times[0,T]$.
This leads to $\partial_t u_m\equiv0$ on $V\setminus\Omega^\circ_m\times[0,T]$. Hence $\partial_t u_m\in\mathcal{H}_m$, and so $u_m-\partial_t u_m\in\mathcal{H}_m$ on $[0,T]$.

In \eqref{e:um_on_Omega}, putting $v=u_m-\partial_t u_m$, we obtain
\begin{eqnarray*}
\|\partial_t u_m(\cdot,t)\|^2_{L^2(\Omega^\circ)}\le\|\Delta u_m(\cdot,t)+f(\cdot,t)\|^2_{L^2(\Omega^\circ)}.
\end{eqnarray*}
Combining this with
\begin{eqnarray*}
\begin{aligned}
\|\Delta u_m(\cdot,t)\|^2_{L^2(\Omega^\circ)}
&\le \sum_{x\in\Omega}\frac{1}{\mu(x)}\Big(\sum_{y\sim x}\omega_{xy}\Big)\Big(\sum_{y\sim x}\omega_{xy}\big|u_m(y,t)-u_m(x,t)\big|^2\Big)\\
&\le 2 D_\mu\sum_{x\in\Omega^\circ}\sum_{y\sim x}\omega_{xy}\big(|u_m(y,t)|^2+|u_m(x,t)|^2\big)\\
&\le 2 D_\mu M_d\sum_{x\in\Omega^\circ}\frac{1}{\mu(x)}\sum_{y\sim x}\omega_{xy}|u_m(x,t)|^2\mu(x)\\
&\le 2 (D_\mu M_d)^2\|u_m(\cdot,t)\|^2_{L^2(\Omega)},
\end{aligned}
\end{eqnarray*}
we get $\partial_t u_m$ is uniformly bounded in $\mathcal{L}^2$. Furthermore, $\partial_t u_m$ is uniformly bounded in $L^2([0,T]; \mathcal{L}^2)$. This completes the proof.
\end{proof}

\begin{lem}\label{L:conve_L_2_var}
There exist a subsequence of $\{u_{m}\}$, denoted also by $\{u_m\}$, and a function $u\in\mathcal{L}^2$ such that for any $t\in[0,T]$,
$$u_{m}(\cdot,t)\to u(\cdot,t),\qquad \Delta u_{m}(\cdot,t)\to \Delta u(\cdot,t), \qquad \partial_t u_{m}(\cdot,t)\to \partial_t u(\cdot,t)\quad\mbox{ in }\mathcal{L}^2.$$
\end{lem}
\begin{proof}
This proof consists three parts.

{\bf Part 1.} $u_{m}\to u$ in $\mathcal{L}^2$.

Lemma \ref{L:uniform_bdd_seq_var} implies that $\|u_{m}(\cdot,t)\|_{L^2(\Omega^\circ)}$ is uniformly bounded on $[0,T]$.
Hence $\|u_{m}(\cdot,t)\|_{L^2(\Omega_1^\circ)}$ is uniformly bounded, and so there exists a subsequence $\{u_{1m}\}\subseteq\{u_{m}\}$ and a function $u^\star_1$ such that $u_{1m}\to u^\star_1$ in $L^2(\Omega_1^\circ)$. It is obvious that $\|u_{1m}(\cdot,t)\|_{L^2(\Omega_1^\circ)}$ is uniformly bounded. Similarly, we can extract a subsequence $\{u_{2m}\}\subseteq\{u_{1m}\}$ such that $u_{2m}\to u^\star_2$ in $L^2(\Omega^\circ_2)$. It is easy to see that $u^\star_2|_{\Omega^\circ_1}=u^\star_1$. Continuing, for
$m=3,\ldots$, we get a subsequence $\{u_{im}\}\subseteq \{u_{i-1 m}\}\subseteq\cdots\subseteq \{u_{2m}\}\subseteq \{u_{1m}\}\subseteq \{u_{m}\}$ and a function $u^\star_i$ such that $u_{im}\to u^\star_i$ in $L^2(\Omega_i^\circ)$. Also, $u^\star_i|_{\Omega_{i-1}^\circ}= u^\star_{i-1}$. Using diagonal method, we can get a subsequence $\{u_{mm}\}$, denote it by $\{u_{m}\}$, and a function $u^\star\in L^2(\Omega^\circ)$ satisfying $u^\star|_{\Omega_m}=u^\star_m$, such that $u_{m}\to u^\star$ in $L^2(\Omega_k^\circ)$ for any bounded $\Omega_k\subseteq \Omega$. Now, we show that $u_{m}\to u^\star$ in $L^2(\Omega^\circ)$. Note that
\begin{eqnarray*}
\lim_{m\to\infty}\int_{\Omega^\circ}|u_m-u^\star|^2\,d\mu
=\lim_{m\to\infty}\lim_{k\to\infty}\int_{\Omega^\circ_k}|u_m-u^\star|^2\,d\mu
=\lim_{k\to\infty}\lim_{m\to\infty}\int_{\Omega^\circ_k}|u_m-u^\star|^2\,d\mu
=0,
\end{eqnarray*}
where we use the facts that $\int_{\Omega^\circ_k}|u_m-u^\star|^2\,d\mu$ is monotonically increasing with respect to $k$, and
$$\int_{\Omega^\circ_k}|u_m-u^\star|^2\,d\mu\le 2\int_{\Omega^\circ_k}(|u_m|^2\,d\mu+|u^\star|^2)\,d\mu\le C+2\|u^\star\|^2_{L^2(\Omega^\circ)}.$$
Let $u:=u^\star$ on $\Omega^\circ$ and $u=0$ on $V\setminus\Omega^\circ$. Then $u_{m}\to u$ in $\mathcal{L}^2$.

{\bf Part 2.} $\Delta u_{m}\to \Delta u$ in $\mathcal{L}^2$.

It is obvious that the result follows from $u_{m}\to u$ in $\mathcal{L}^2$ and
\begin{eqnarray*}
\|\Delta (u_m-u)(\cdot,t)\|^2_{L^2(\Omega^\circ)}\le 2(D_\mu M_d)^2\|(u_m-u)(\cdot,t)\|^2_{L^2(\Omega^\circ)}.
\end{eqnarray*}

{\bf Part 3.} $\partial_t u_{m}\to \partial_t u$ in $\mathcal{L}^2$.

As same as Part 1, we can extract a subsequence of $\{\partial_t u_m\}$, denoted also by $\{\partial_t u_m\}$, such that
$\partial_t u_{m}(\cdot,t)\to S^\ast(\cdot,t)$ in $\mathcal{L}^2$ for some function $S^\ast\in\mathcal{L}^2$. In the following, we need to show that $S^\ast=\partial_t u$. There are two steps.

{\bf Step 1.} We prove that $\partial_t u$ exists. The proof is same to that of Lemma \ref{L:point_conve_eq}(b).
It follows from Lemma \ref{L:uniform_bdd_seq_var} that we can extract a subsequence of $\{\partial_t u_{m}\}$, denoted also by itself, such that for some function $S\in L^2([0,T]; \mathcal{L}^2)$,
\begin{eqnarray*}
\int_0^T(\partial_t u_{m}-S,v)\,dt\to0\qquad\mbox{ for any }v\in L^2([0,T];\mathcal{L}^2).
\end{eqnarray*}
Note that for any $v\in\mathcal{L}^2$,
\begin{eqnarray}\label{e:partial_s_u_m_ell_var}
(u_{m}(\cdot,T)-g_{m},v)=\int_0^T\big(\partial_t u_{m}(\cdot,t),v\big)\,dt.
\end{eqnarray}
Taking limit as $m\to\infty$ in \eqref{e:partial_s_u_m_ell_var} and using $u_{m}\to u$ in $\mathcal{L}^2$, we obtain
\begin{eqnarray*}
\begin{aligned}
(u(\cdot,T)-g,v)=\int_0^T(S(\cdot,t),v)\,dt=\Big(\int_0^TS(\cdot,t)\,dt,v\Big),
\end{aligned}
\end{eqnarray*}
where we use facts that $\lim_{m\to\infty}(g_{m},v)=(g,v)$, and
\begin{equation*}\label{e:fact_2_var}
\int_0^T\Big(\sum_{x\in\Omega^\circ}S(x,t)v(x,t)\mu(x)\Big)\,dt=\sum_{x\in\Omega^\circ}\Big(\int_0^TS(x,t)\,dt\Big)v(x,t)\mu(x).
\end{equation*}
Hence $S$ is integrate on $[0,T]$, $u$ is continuous with respect to $t$. Moreover, $u(x,t)=g(x)$ on $\Omega^\circ$, and $S=\partial_tu$.

{\bf Step 2.} We show that $S^\ast=\partial_t u$.
The fact $u_{m}\to u$ in $\mathcal{L}^2$ implies that $u_{m}\to u$ on $\Omega^\circ\times[0,T]$.
Hence, for any $t_0\in[0,T]$ and $x\in\Omega^\circ$,
\begin{eqnarray*}
\lim_{m\to\infty}\partial_t u_{m}(x,t_0)=\partial_t u(x,t_0).
\end{eqnarray*}
That is, $\partial_t u_{m}\to \partial_t u$ on $\Omega^\circ\times [0,T]$.

On the other hand, since $\partial_t u_{m}\to S^\ast$ in $\mathcal{L}^2$, we get
$\partial_t u_{m}\to S^\ast$ on $\Omega^\circ\times[0,T]$. Hence $S^\ast=\partial_t u$ on $\Omega^\circ\times [0,T]$.
This completes the proof.
\end{proof}

Using Lemma \ref{L:conve_L_2_var}, we can get the following results.
\begin{lem}\label{L:conve_L_2_int_var}
There exist a subsequence of $\{u_m\}$ and $u\in\mathcal{L}^2$ such that for any $t\in[0,T]$,
\begin{enumerate}
\item[(a)] for any $v\in \mathcal{L}^2$,
$$
\int_{\Omega^\circ}\partial_t u_{m}\cdot v\,d\mu\to \int_{\Omega^\circ}\partial_t u\cdot v\,d\mu\quad\mbox{ and }\quad
\int_{\Omega^\circ}\Delta u_{m}\cdot v\,d\mu\to \int_{\Omega^\circ}\Delta u\cdot v\,d\mu;
$$
\item[(b)] $\int_{\Omega^\circ}\partial_t u_{m}\cdot u_{m}\,d\mu\to \int_{\Omega^\circ}\partial_t u\cdot u\,d\mu$,
$\int_{\Omega^\circ}\Delta u_{m}\cdot u_{m}\,d\mu\to \int_{\Omega^\circ}\Delta u\cdot u\,d\mu$,
\item[(c)] $\int_\Omega f_{m}\cdot u_{m_j}\,d\mu\to \int_\Omega f\cdot u\,d\mu$.
\end{enumerate}
\end{lem}

\begin{proof}
(a) Using Lemma \ref{L:conve_L_2_var} and H\"older inequality, it is easy to see that for any $v\in\mathcal{L}^2$,
\begin{eqnarray*}
\int_{\Omega^\circ} (\partial_t u_{m}-\partial_t u)\cdot v\,d\mu
\le\Big(\int_{\Omega^\circ} |\partial_t u_{m}-\partial_t u|^2\,d\mu\Big)^{1/2}\cdot \Big(\int_{\Omega^\circ} |v|^2\,d\mu\Big)^{1/2}\to 0.
\end{eqnarray*}
Hence
$\int_{\Omega^\circ}\partial_t u_{m}\cdot v\,d\mu\to \int_{\Omega^\circ}\partial_t u\cdot v\,d\mu$. Similarly, we get
$\int_{\Omega^\circ}\Delta u_{m}\cdot v\,d\mu\to \int_{\Omega^\circ}\Delta u\cdot v\,d\mu$.

(b) follows from Lemma \ref{L:conve_L_2_var}, triangle inequality and H\"older inequality.

(c) It follows from \eqref{e:condition_for_f} that $f\in L^2(\Omega^\circ)$ on $[0,\infty)$. Combining this with
$f_m\to f$ on $\Omega^\circ\times[0,\infty)$, using Dominated Convergence Theorem, we obtain $f_m\to f$ in $L^2(\Omega^\circ)$.
Combining this with triangle and H\"older inequalities, we get the result holds.
\end{proof}

In the rest of this section, using the notation defined above, we give the proof of Theorem~\ref{T:var_ineq}.
\begin{proof}[Proof of Theorem \ref{T:var_ineq}]
Lemma \ref{L:unique_solution_on_Omega_m} implies that for any $T>0$, any $t\in[0,T]$ and any $v\in\mathcal{H}_m$,
\begin{eqnarray*}
\int_{\Omega^\circ} \partial_t u_m\cdot(v-u_m)\,d\mu\ge \int_{\Omega^\circ}(\Delta u_m+f_m)\cdot(v-u_m)\,d\mu,
\end{eqnarray*}
and
\begin{eqnarray*}
u_m(x,0)=g_m(x)\quad\mbox{ for }x\in \Omega_m^\circ,
\qquad\mbox{and}\qquad u_m(x,t)=0\quad\mbox{ for }(x,t)\in V\setminus\Omega_m^\circ\times[0,T].
\end{eqnarray*}
Using Lemma \ref{L:conve_L_2_int_var}, the proof is completed.
\end{proof}

\textbf{Acknowledgement}
This research was supported by the National Science Foundation of China [grants 12071245].

Yong Lin,\\
Yau Mathematical Sciences Center, Tsinghua University\\ Beijing, 100084, China.\\
\textsf{E-mail: yonglin@tsinghua.edu.cn}\\
\\
Yuanyuan Xie,\\
School of Mathematics, Renmin University of China\\ Beijing, 100872, China.\\
\textsf{E-mail: yyxiemath@163.com}\\

\end{document}